\documentclass{article}
\usepackage{amsmath,amssymb,amsfonts,amsthm,amsopn,dsfont}
\usepackage{bm}
\usepackage{subcaption}
\usepackage{graphicx}
\usepackage{float}

\usepackage[utf8]{inputenc}
\usepackage[T1]{fontenc}

\usepackage{xcolor}

\usepackage[ruled,vlined,linesnumbered]{algorithm2e}
\usepackage{multirow}

\newcommand{\hsp}{\mathcal{H}^{\Gamma}}

\newcommand{\extg}{\mathcal{E}_{_{\Gamma}}}
\newcommand{\exts}{\mathcal{E}_{_{\Sigma}}}
\newcommand{\trg}{\gamma_{_{\Gamma}}}

\begin{document}


	\title{A gradient based resolution strategy for a PDE-constrained optimization approach for 3D-1D coupled problems}

	\author{Stefano Berrone
	\footnote{stefano.berrone@polito.it}, 
	Denise Grappein
	\footnote{denise.grappein@polito.it},
	Stefano Scial\'o
    \footnote{stefano.scialo@polito.it},
	Fabio Vicini
	\footnote{fabio.vicini@polito.it}\\{\footnotesize Dipartimento di Scienze Matematiche, Politecnico di Torino,
	Corso Duca degli Abruzzi 24,}\\{\footnotesize 10129 Torino, Italy. Member of INdAM research group GNCS.}
}
\date{}

	
\maketitle
	\section*{Abstract}
		Coupled 3D-1D problems arise in many practical applications, in an attempt to reduce the computational burden in simulations where cylindrical inclusions with a small section are embedded in a much larger domain. Nonetheless the resolution of such problems can be non trivial, both from a mathematical and a geometrical standpoint. Indeed 3D-1D coupling requires to operate in non standard function spaces, and, also, simulation geometries can be complex for the presence of multiple intersecting domains. Recently, a PDE-constrained optimization based formulation has been proposed for such problems, proving a well posed mathematical formulation and allowing for the use of non conforming meshes for the discrete problem. Here an unconstrained optimization formulation of the problem is derived and an efficient gradient based solver is proposed for such formulation. Some numerical tests on quite complex configurations are discussed to show the viability of the method.


	\subsection*{Keywords}
		3D-1D coupling - three-field - domain-decomposition - non conforming mesh - optimization methods for elliptic problems




\section{Introduction}
This work presents a conjugate gradient based resolution strategy for a recently developed numerical scheme for the coupling of three-dimensional and one-dimensional elliptic equations (3D-1D coupling) \cite{3D1Darxiv}. Coupled problems with such dimensionality gap arise, in particular, when small tubular inclusions  embedded in a much wider domain are dimensionally reduced to 1D manifolds for computational efficiency. This allows to avoid the complexity related to the building of a three-dimensional grid within the inclusions. Examples of applications range from the description of biological tissues \cite{Notaro2016,Koppl2020}, roots-soil interaction \cite{Schroder2012,Koch2018}, geological reservoir simulations \cite{Gjerde2018a,Gjerde2020, CerLauZun2019}, to fiber-reinforced materials \cite{Steinbrecher2020,Llau2016}.

The mathematical treatment of the coupling between a 3D and a 1D problem is non trivial, as no bounded trace operator is defined when the dimensionality gap between the interested manifolds is higher than one. In \cite{Dangelo2012} suitable weighed Sobolev spaces were introduced, thanks to which a bounded trace operator was defined and the well-posedness of the problem was worked out by means of the Banach-Ne\v{c}as-Babu\v{s}ka theorem \cite{ErnGuermond}. Other approaches rely on the use of regularizing techniques \cite{Tornberg2004} or lifting strategies \cite{Koppl2018}. In \cite{Zunino2019} a topological model reduction is employed and averaging operators are introduced leading to a well posed 3D-1D coupled problem. Problems with singular sources on lines are also studied in \cite{Gjerde2019}, where an approach based on the splitting of the solution in a low regularity part and a regular correction is analysed.

The present work is based on a re-formulation of the original 3D-3D problem into properly defined functional spaces, thus paving the way for a well posed formulation of the reduced 3D-1D problem \cite{3D1Darxiv}. The numerical resolution is further obtained through a PDE-constrained optimization based approach  \cite{BPSc,BSV,BDS,ThreeFieldSISC}, in which problems in the 3D bulk domain and in the 1D inclusions are decoupled using a three-field based domain decomposition method.
A cost functional, expressing the error in the fulfilment of interface conditions, is minimized to restore the coupling.
The discrete problem is re-written as an unconstrained optimization problem and a conjugate gradient scheme is proposed for its numerical resolution. This allows to treat efficiently large scale problems.

The manuscript is organized as follows: the problem of interest is briefly recalled in Section~\ref{not_and_form}, along with its re-formulation as a PDE-constrained optimization problem. The corresponding discrete version is described in Section~\ref{Discrete}, whereas the novel conjugate gradient based resolution strategy is presented in Section~\ref{sec:CG}. Three numerical tests are provided in Section~\ref{Num_res}, and some conclusions are proposed in Section~\ref{sec:Concl}.

\section{Notation and problem formulation}\label{not_and_form}
\begin{figure}
	\centering
	\includegraphics[width=0.55\textwidth]{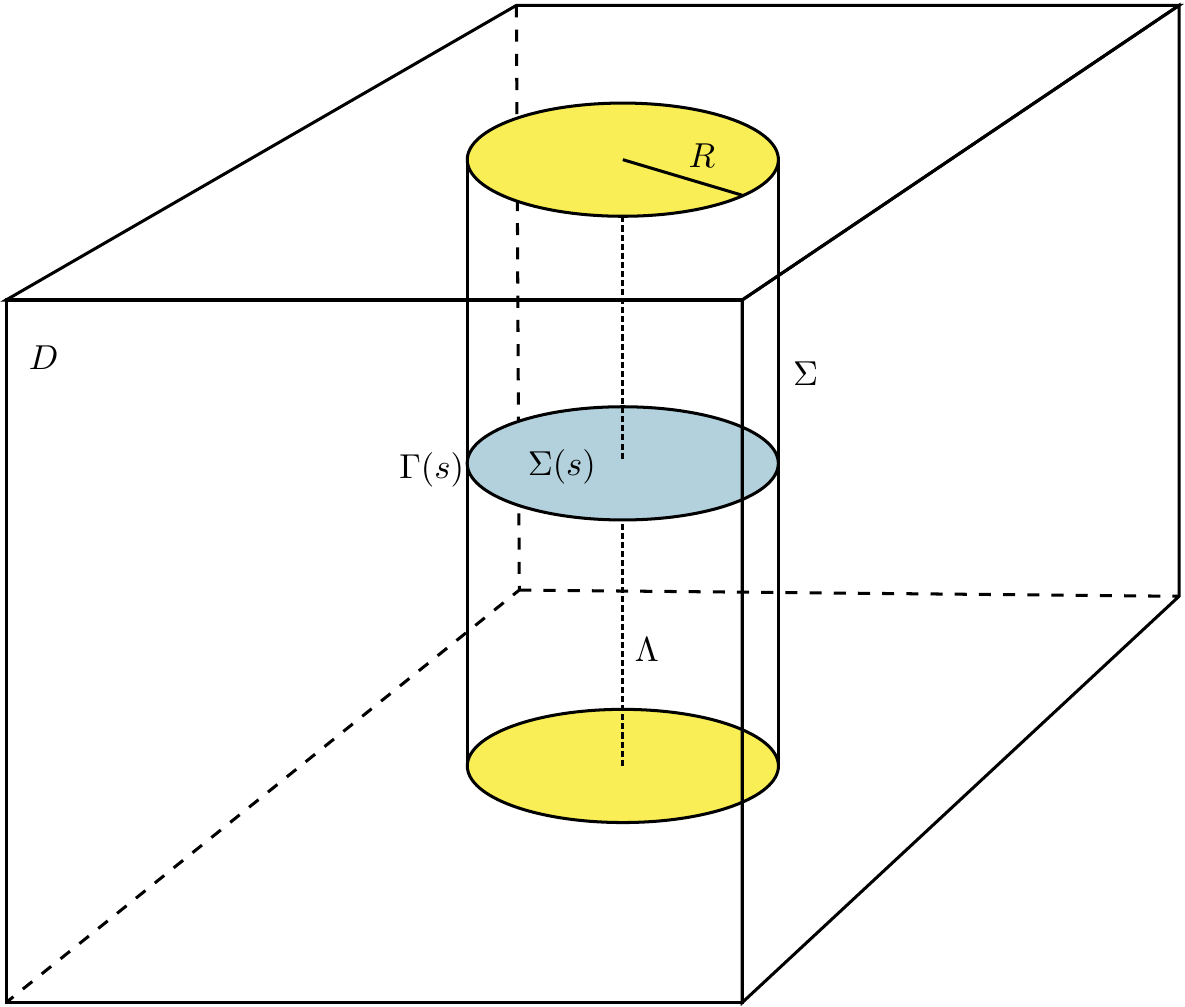}
	\caption{Example of Domain $\Omega$}
	\label{fig:domain}
\end{figure}
We briefly recall here the derivation of the reduced 3D-1D coupled problem from the original equi-dimensional formulation, referring to \cite{3D1Darxiv} for a more comprehensive discussion.

Let us consider a three dimensional convex domain $\Omega$ with a single cylindrical inclusion $\Sigma \in \mathbb{R}^3$ with centreline $\Lambda=\left\lbrace\bm{\lambda}(s), s \in (0,S) \right\rbrace$, see Figure~\ref{fig:domain}. We denote by $\Sigma(s)$ the transverse section of $\Sigma$ at $s\in [0,S]$ of radius $R \ll \text{diam}(\Omega)$ and by $\Gamma(s)$ its boundary. The lateral surface of the whole cylinder is $\Gamma$, whereas $\Sigma_0=\Sigma(0)$ and $\Sigma_S=\Sigma(S)$ are the two extreme sections. Let us set $D=\Omega \setminus \Sigma$ the domain without the cylindrical inclusion and let us denote by $\partial D=\partial \Omega \cup \left\lbrace \Gamma\cup \Sigma_0 \cup \Sigma_S\right\rbrace $ its boundary, being $\partial \Omega$ the boundary of $\Omega$. For simplicity of exposition we assume here that $\Sigma_0$ and $\Sigma_S$ lie on $\partial \Omega$, and thus we introduce the symbol $\partial D^e=\partial \Omega \setminus \left\lbrace \Sigma_0 \cup \Sigma_S\right\rbrace$ to denote the external boundary of domain $D$.
We are interested in the following problem in $\Omega$:
\begin{align}
	- & \nabla \cdot (\bm{K} \nabla u)=f                         & \text{in } D\label{eqOmega}                        \\[0.35em]
	- & \nabla \cdot(\tilde{\bm{K}}\nabla\tilde{u})=g            & ~\text{ in } \Sigma \label{eqSigma}                \\
	  & u=0                                                      & \text{on } \partial D^e \label{cbOmega}            \\
	  & u_{|_{\Gamma}}=\psi                                      & ~\text{ on } \Gamma \label{psi_u}                  \\
	  & \bm{K}\nabla u \cdot \bm{n}=\phi                         & ~\text{on } \Gamma \label{phi_u}                   \\
	  & \tilde{u}=0                                              & \text{ on } \Sigma_0 \cup \Sigma_S \label{cbSigma} \\
	  & \tilde{u}_{|_{\Gamma}}=\psi                              & ~\text{ on } \Gamma \label{psi_utilde}             \\
	  & \tilde{\bm{K}}\nabla \tilde{u}\cdot \bm{\tilde{n}}=-\phi & ~\text{on } \Gamma\label{phi_utilde}
\end{align}
where $u$ and $\tilde{u}$ are the unknowns related to domains $D$ and $\Sigma$, respectively, $\bm{n}$ and $\bm{\tilde{n}}$ are unit normal vectors to $\Gamma$ outward pointing from $D$ and $\Sigma$, respectively, $\bm{K}$ and $\tilde{\bm{K}}$ are uniformly positive definite tensors in $D$ and $\Sigma$, respectively, and $f$ and $g$ are source terms. Equations \eqref{psi_u},\eqref{psi_utilde} and \eqref{phi_u},\eqref{phi_utilde}, namely the pressure continuity and the flux conservation conditions on the interface $\Gamma$, could be written as $u_{|_{\Gamma}}=\tilde{u}_{|_{\Gamma}}$ and $\bm{K}\nabla u =-\tilde{\bm{K}}\nabla \tilde{u}\cdot \bm{\tilde{n}} $. Nevertheless the equations can be split, as shown above, by introducing the auxiliary variables $\phi$ and $\psi$, in view of the application of a three-field domain decomposition approach.

As mentioned, when $R$ is much smaller than the domain size, it can be computationally convenient to recast the previous problem in a 3D-1D coupled problem, assuming that the variations of $\tilde{u}$ on the cross sections of the cylinder can be considered negligible, as well as the variations of $\psi$ on $\Gamma (s)$. 
In order to derive a well posed 3D-1D coupled problem, we introduce the following function spaces:
\begin{align*}
	 & H_0^1(D)=\left\lbrace v \in H^1(D) : v_{|_{\partial D^e}}=0 \right\rbrace,                      \\
	 & H_0^1(\Sigma)=\left\lbrace v \in H^1(\Sigma): v_{|_{\Sigma_0}}=v_{|_{\Sigma_S}}=0\right\rbrace, \\
	 & H_0^1(\Lambda)=\left\lbrace v \in H^1(\Lambda): v(0)=v(S)=0\right\rbrace,
\end{align*}
the trace operator $\gamma_{_\Gamma}:H^1(D)\cup H^1(\Sigma)\rightarrow H^{\frac{1}{2}}(\Gamma)$ s.t.
\begin{equation}
	\gamma_{_\Gamma}v=v_{|_\Gamma} ~\forall v \in H^1(D)\cup H^1(\Sigma)
\end{equation}
and the two extension operators $\mathcal{E}_{_\Sigma}: H^1(\Lambda) \rightarrow H^1(\Sigma) $ and $\mathcal{E}_{_\Gamma}: H^1(\Lambda) \rightarrow H^{\frac{1}{2}}(\Gamma)$ such that for any $\hat{v} \in H_0^1(\Lambda)$ $\mathcal{E}_{_\Sigma}(\hat{v})$ is the extension of the point-wise value $\hat{v}(s)$, $s \in [0,S]$, to $\Sigma(s)$ and $\mathcal{E}_{_\Gamma}(\hat{v})$ is the extension of $\hat{v}(s)$ to $\Gamma(s)$. Let us observe that $\mathcal{E}_{_\Gamma}=\gamma_{_\Gamma}\circ \mathcal{E}_{_\Sigma}$.
Let us further consider the spaces:
\begin{align*}
	 & \hat{V}=H_0^1(\Lambda),                                                                                    \\
	 & \widetilde{V}=\lbrace v \in H_0^1(\Sigma): v =\exts\hat{v}, ~\hat{v} \in \hat{V} \rbrace,                  \\
	 & \mathcal{H}^{\Gamma}=\lbrace v \in H^{\frac{1}{2}}(\Gamma): v =\extg\hat{v}, ~\hat{v} \in \hat{V} \rbrace, \\
	 & V_D=\left\lbrace v \in H_0^1(D): \gamma_{_\Gamma}v \in \mathcal{H}^{\Gamma}\right\rbrace.
\end{align*}
We can observe that functions in $\widetilde{V}$ are the extension to the whole domain $\Sigma$ of functions defined on the centreline $\Lambda$.
Similarly functions in $\mathcal{H}^{\Gamma}$ are extension on $\Gamma$ of functions in $\hat{V}$ or, equivalently, traces on $\Gamma$ of element in $\widetilde{V}$. 
The space $V_D$ contains functions whose trace on $\Gamma$ belongs to $\mathcal{H}^{\Gamma}$.

Denoting by $(\cdot,\cdot)_{\star}$ the $L^2$-scalar product on a generic domain $\star$ and indicating with $X'$ the dual of a generic space $X$, we can write a well posed weak formulation of problem  \eqref{eqOmega}-\eqref{phi_utilde} in the above function spaces as follows:
\textit{find} $(u, \tilde{u}) \in V_D \times \widetilde{V}$, $\phi\in {\mathcal{H}^{\Gamma}}'$ and $\psi \in \mathcal{H}^{\Gamma}$ \textit{such that:}
\begin{align}
	 & (\bm{K}\nabla u, \nabla v)_{D}-\left\langle \phi,\gamma_{_\Gamma}v \right\rangle_{{\mathcal{H}^{\Gamma}}', {\mathcal{H}^{\Gamma}}}=(f,v)_{D}~                                                    & \forall v \in V_D,~\phi \in {\mathcal{H}^{\Gamma}}' \label{eq_var_u}                      \\
	 & (\bm{\tilde{K}}\nabla \tilde{u}, \nabla \tilde{v})_{{\Sigma}}+\left\langle \phi,\gamma_{_\Gamma}\tilde{v} \right\rangle_{{\mathcal{H}^{\Gamma}}', {\mathcal{H}^{\Gamma}}}=(g,\tilde{v})_{\Sigma} & \forall \tilde{v} \in \widetilde{V},~\phi \in {\mathcal{H}^{\Gamma}}' \label{eq_var_uhat}
\end{align}
\begin{align}
	 & \left\langle \gamma_{_\Gamma}u-\psi,\eta\right\rangle_{\mathcal{H}^{\Gamma},{\mathcal{H}^{\Gamma}}'}=0         & ~\forall \eta \in  {\mathcal{H}^{\Gamma}}', \psi \in \mathcal{H}^{\Gamma} \label{condpsi_u}  \\
	 & \left\langle \gamma_{_\Gamma}\tilde{u}-\psi,\eta\right\rangle_{\mathcal{H}^{\Gamma},{\mathcal{H}^{\Gamma}}'}=0 & ~\forall \eta \in  {\mathcal{H}^{\Gamma}}', \psi \in \mathcal{H}^{\Gamma}\label{condpsi_hat}
\end{align}
where $\phi$ is the unknown flux through $\Gamma$ and $\psi$ represents the value of the solution on $\Gamma$.
We are now interested in solving this problem, which has the advantage that it can be easily recast in a 3D-1D reduced problem while still working with a well posed trace operator  $\gamma_{_\Gamma}(\cdot)$ from a three-dimensional to a two dimensional manifold.
Recalling that:
\begin{equation*}
	\left\langle \phi,\gamma_{_\Gamma}v \right\rangle_{{\mathcal{H}^{\Gamma}}', {\mathcal{H}^{\Gamma}}}=\int_0^S\Big( \int_{\Gamma(s)}\phi~\gamma_{_\Gamma}v~dl\Big) ds \quad \forall v \in V_D
\end{equation*}
we have, denoting by $|\Gamma(s)|$ the perimeter of the section at $s\in[0,S]$ and by $\overline{\phi}(s)$ the mean value of $\phi$ on $\Gamma(s)$, that
\begin{equation*}
	\left\langle \phi,\gamma_{_\Gamma}v \right\rangle_{{\mathcal{H}^{\Gamma}}', {\mathcal{H}^{\Gamma}}}=\int_0^S|\Gamma(s)|\overline{\phi}(s)\check{v}(s)~ds=\left\langle |\Gamma|\overline{\phi},\check{v}\right\rangle_{\hat{V}',\hat{V}}.
\end{equation*}
Function $\check{v} \in \hat{V}$ is introduced s.t. given $v \in V_D$ for all $s \in [0,S]$ we have by definition $\gamma_{_\Gamma}v=\extg\check{v}=\check{v}(s)$; thus, $\int_{\Gamma(s)}\phi\ \trg v  ~dl= \check{v}(s) \int_{\Gamma(s)}\phi ~dl= \check{v}(s) |\Gamma(s)| \overline{\phi}(s)$.

Proceeding in a similar way we can rewrite equations \eqref{condpsi_u} and \eqref{condpsi_hat} as
\begin{align*}
	 & \left\langle \gamma_{_\Gamma}u-\psi,\eta\right\rangle_{\mathcal{H}^{\Gamma},{\mathcal{H}^{\Gamma}}'}=\left\langle |\Gamma|(\check{u}-\hat{\psi}),\overline{\eta}\right\rangle_{\hat{V},\hat{V}'}=0      \\
	 & \left\langle \gamma_{_\Gamma}\tilde u-\psi,\eta\right\rangle_{\mathcal{H}^{\Gamma},{\mathcal{H}^{\Gamma}}'}=\left\langle |\Gamma|(\hat{u}-\hat{\psi}),\overline{\eta}\right\rangle_{\hat{V},\hat{V}'}=0
\end{align*} where $\check{u},\hat{\psi}\in \hat{V}$ are such that $\gamma_{_\Gamma}u=\extg\check{u}$, $\psi=\extg\hat{\psi}$ and $\gamma_{_\Gamma}\tilde{u}=\gamma_{_\Gamma}\exts\tilde{u}=\extg\hat{u}$, as $\tilde{u} \in \widetilde{V}$.
Concerning the problem in $\Sigma$:
\begin{equation*}
	(\bm{\tilde{K}}\nabla \tilde{u}, \nabla \tilde{v})_{\Sigma}=\int_{\Sigma}\bm{\tilde{K}}\nabla \tilde{u}\nabla\tilde{v}~d\sigma=\int_0^S\bm{\tilde{K}}|\Sigma(s)|\cfrac{d\hat{u}}{ds}~\cfrac{d\hat{v}}{ds}~ds
\end{equation*}
where $\hat{u},\hat{v} \in \hat{V}$ are such that $\tilde{u}=\exts\hat{u}$, $\tilde{v}=\exts\hat{v}$ and $|\Sigma(s)|$ is the section area at $s \in [0,S]$.

If we denote by $V$ the space obtained extending $V_D$ from $D$ to the whole region $\Omega$, we can set the limit problem (\ref{eq_var_u})-(\ref{condpsi_hat}) as a reduced 3D-1D coupled problem:
\textit{Find $(u,\hat{u}) \in V\times\hat{V}$, $\overline{\phi}\in \hat{V}'$ and $\hat{\psi} \in \hat{V}$ such that:}
\begin{align}
	 & (\bm{{K}}\nabla u, \nabla v)_{\Omega}-\left\langle |\Gamma|\overline{\phi},\check{v} \right\rangle_{\hat{V}',\hat{V}}=(f,v)_{\Omega}  \quad~\forall v \in V, \check{v} \in \hat{V}: \gamma_{_\Gamma} v=\mathcal{E}_{_\Gamma}\check{v} \label{equaz1}                 \\
	 & \Big( \bm{\tilde{K}}|\Sigma|\frac{d\hat{u}}{ds},\frac{d\hat{v}}{ds}\Big)_{\Lambda}+\left\langle |\Gamma|\overline{\phi},\hat{v} \right\rangle_{\hat{V}',\hat{V}}=(|\Sigma|\overline{\overline{g}},\hat{v})_{\Lambda}\quad~ \forall \hat{v} \in \hat{V}\label{equaz2} \\[0.5em]
	 & \qquad\left\langle |\Gamma|(\check{u}-\hat{\psi}),\overline{\eta}\right\rangle_{\hat{V}',\hat{V}}=0 \qquad \gamma_{_\Gamma}u=\mathcal{E}_{\Gamma}\check{u},~\forall \overline{\eta} \in \hat{V}'\label{condiz1}                                                      \\
	 & \qquad\left\langle |\Gamma|(\hat{u}-\hat{\psi}),\overline{\eta}\right\rangle_{\hat{V}',\hat{V}}=0 \qquad\forall \overline{\eta} \in \hat{V}'\label{condiz2}
\end{align}
with $\overline{\overline{g}}(s)=\frac{1}{|\Sigma(s)|}\int_{\Sigma(s)}g~d\sigma$, for a sufficiently regular $g$.

Problem \eqref{equaz1}-\eqref{condiz2} can be conveniently stated as a PDE-constrained optimization problem, which yields a discrete problem that can be efficiently solved on independent meshes for the 3D and 1D domains through a gradient based iterative solver. At this end, let us introduce the functional
\begin{eqnarray}
	J(\overline{\phi},\hat{ \psi})&=&\cfrac{1}{2}\left( ||\gamma_{_\Gamma}u(\overline{\phi},\hat{ \psi})-\psi||_{\hsp}^2+||\gamma_{_\Gamma}\tilde{u}(\overline{\phi},\hat{ \psi})-\psi||_{\hsp}^2\right) \nonumber \\
	&=&\cfrac{1}{2}\left( ||\gamma_{_\Gamma}u(\overline{\phi},\hat{ \psi})-\extg \hat{\psi}||_{\hsp}^2+||\gamma_{_\Gamma}\exts \hat{u}(\overline{\phi},\hat{ \psi})-\extg \hat{\psi}||_{\hsp}^2\right)	\label{functional}
\end{eqnarray}
expressing the error in the fulfilment of conditions \eqref{condiz1}-\eqref{condiz2}. Equations \eqref{equaz1}-\eqref{equaz2} are slightly modified as follows:
\begin{align}
	 & (\bm{{K}}\nabla u, \nabla v)_{\Omega}+\alpha(|\Gamma|\check{u},\check{v} )_{\Lambda}-\left\langle |\Gamma|\overline{\phi},\check{v} \right\rangle_{\hat{V}',\hat{V}} =(f,v)_{\Omega} + \alpha(|\Gamma|\hat{\psi},\check{v} )_{\Lambda} \label{eq_stab_u}                                                                                    \\[-0.4em]
	 & \hspace{7cm}\forall v \in V, ~\check{v} \in \hat{V}: \trg v=\extg\check{v}, \nonumber                                                                                                                                                                                                                                                       \\
	 & \Big(\bm{\tilde{K}}|\Sigma|\cfrac{d\hat{u}}{ds},\cfrac{d\hat{v}}{ds}\Big)_{\Lambda}+\hat{\alpha}(|\Gamma|\hat{u},\hat{v})_{\Lambda}+\left\langle |\Gamma|\overline{\phi},\hat{v} \right\rangle_{\hat{V}',\hat{V}}=(|\Sigma|\overline{\overline{g}},\hat{v})_{\Lambda}+\hat{\alpha}(|\Gamma|\hat{\psi},\hat{v})_{\Lambda}\label{eq_stab_hat} \\[-0.4em]
	 & \hspace{8cm}\qquad\forall \hat{v} \in \hat{V}\nonumber.
\end{align}
where the consistent corrections depending from the parameters $\alpha, \hat{\alpha} > 0$ are introduced in order to guarantee the well posedness of the problems independently written on the various domains.

Problem \eqref{equaz1}-\eqref{condiz2} then becomes:
\begin{equation}
	\min_{(\overline{\phi},\hat{\psi})}J(\overline{\phi},\hat{\psi}) \quad \text{ subject to \eqref{eq_stab_u}-\eqref{eq_stab_hat}} \label{minJ}
\end{equation}

\newcommand{\I}{\mathcal{I}}

\section{Matrix formulation}\label{Discrete}
\begin{figure}
	\centering
	\includegraphics[width=0.55\textwidth]{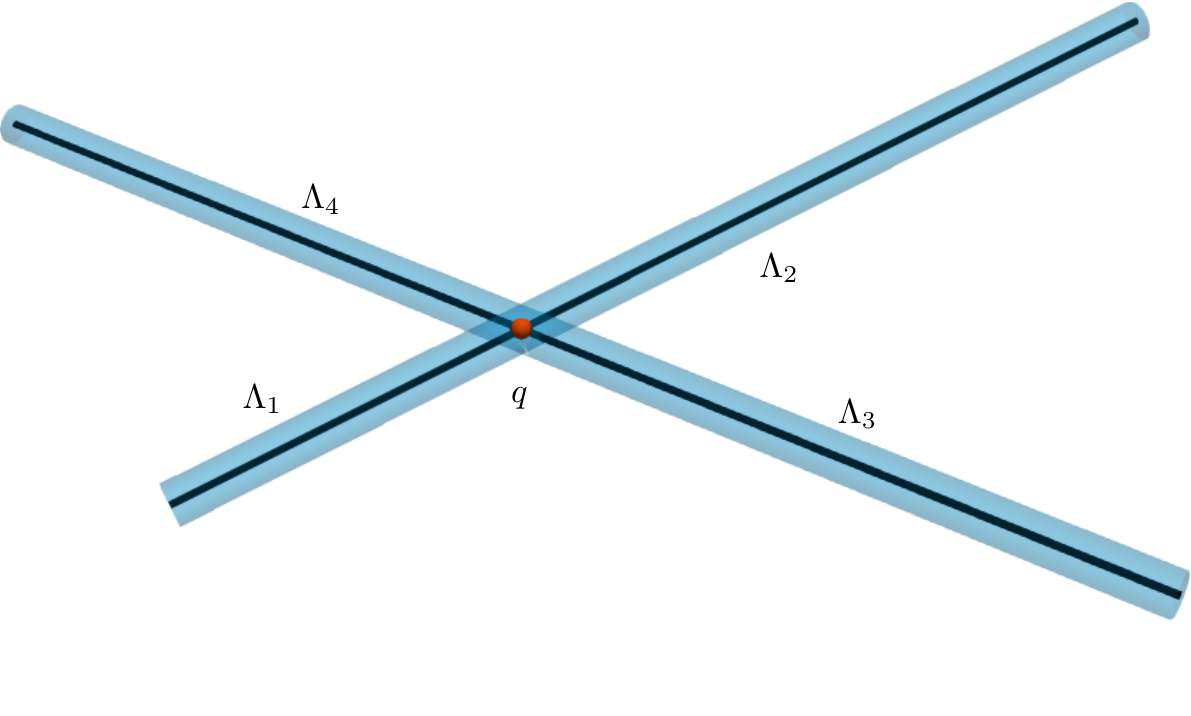}
	\caption{Four segments intersecting at one endpoint $q$}
	\label{seg_inters}
\end{figure}
Let us now derive the discrete counterpart of problem \eqref{minJ}, and thus, let us introduce a tetrahedral mesh $\mathcal{T}$ of domain $\Omega$, and linear Lagrangian finite element basis functions $\left\lbrace \varphi_k \right\rbrace_{k=1}^{N}$ on the mesh $\mathcal{T}$. We will now take into account the more general case where $\I$ segments are embedded in $\Omega$. Intersecting and branching segments are considered as independent segments meeting at one of their end-points, see Figure~\ref{seg_inters}. We build three (possibly) different one-dimensional meshes on each segment $\Lambda_i$, $i=1\ldots,\I$, and we denote them by $\hat{\mathcal{T}_i}$, $\tau^{\phi}_i$ and $\tau^{\psi}_i$. These meshes are independent from each other and from the three-dimensional grid $\mathcal{T}$. We then introduce on such meshes the following basis functions: $\left\lbrace\hat{\varphi}_{i,k} \right\rbrace _{k=1}^{\hat{N}_i}$ on $\hat{\mathcal{T}_i}$, $\left\lbrace \theta_{i,k}\right\rbrace_{k=1}^{N_i^{\phi}}$ on $\tau^{\phi}_i$ and $\left\lbrace \eta_{i,k}\right\rbrace_{k=1}^{N_i^{\psi}}$ on $\tau^{\psi}_i$. We have
\begin{equation*}
	U=\sum_{k=1}^{N}U_k\varphi_k, \quad	\hat{U}_i=\sum_{k=1}^{\hat{N}_i}\hat{U}_{i,k}\hat{\varphi}_{i,k}, \quad \Phi_i=\sum_{k=1}^{N_i^{\phi}}\Phi_{i,k}\theta_{i,k}, \quad \Psi_i=\sum_{k=1}^{N_i^{\psi}}\Psi_{i,k}\eta_{i,k}
\end{equation*}
representing the discrete versions of variable $u$ in $\Omega$, and $\hat{u}_i$, $\overline{\phi}_i$, $\hat{\psi}_i$ on each segment $\Lambda_i$, $i=1,\ldots,\I$.
Further, replacing the definitions of the discrete variables into the constraint equations, we collect the integrals of basis functions into matrices as follows:
\begin{align*}
	 & \bm{A} \in \mathbb{R}^{N\times N} \text{ s.t. } (A)_{kl}=\left(\bm{K}\nabla\varphi_k,\nabla\varphi_l\right)_{\Omega}+\alpha\sum_{i=1}^{\I}\left(|\Gamma(s_i)|{\varphi_k}_{|_{\Lambda_i}},{\varphi_l}_{|_{\Lambda_i}}\right)_{{\Lambda_i}}                                                                 \\[0.5em]
	 & \bm{\hat{A}_i} \in \mathbb{R}^{\hat{N}_i\times \hat{N}_i} \text{ s.t. } (\hat{A}_i)_{kl}=\left(\bm{\tilde{K}_i}|\Sigma(s_i)|\frac{d\hat{\varphi}_{i,k}}{ds},\frac{d\hat{\varphi}_{i,l}}{ds}\right)_{\Lambda_i} +\hat{\alpha}\left(|\Gamma(s_i)|\hat{\varphi}_{i,k},\hat{\varphi}_{i,l}\right)_{\Lambda_i}
\end{align*}
\begin{align*}
	 & \bm{B_i} \in \mathbb{R}^{N\times N_i^{\phi}} \text{ s.t. } (B_i)_{kl}=\left(|\Gamma(s_i)|{{\varphi_{k}}_{|_{\Lambda_i}},\theta_{i,l}}\right)_{\Lambda_i}                                       \\
	 & \bm{\hat{B}_i} \in \mathbb{R}^{\hat{N}_i\times N_i^{\phi}} \text{ s.t. } (\hat{B}_i)_{kl}=\left(|\Gamma(s_i)|{\hat{\varphi}_{i,k},\theta_{i,l}}\right)_{\Lambda_i}                             \\
	 & \bm{C_i}^{\alpha} \in \mathbb{R}^{N\times N_i^{\psi}} \text{ s.t. } (C_i^{\alpha})_{kl}=\alpha\left(|\Gamma(s_i)|{\varphi_k}_{|_{\Lambda_i}},\eta_{i,l}\right)_{\Lambda_i}                     \\
	 & \bm{\hat{C}_i}^{\alpha} \in \mathbb{R}^{\hat{N}_i\times N_i^{\psi}} \text{ s.t. } (\hat{C_i}^{\alpha})_{kl}=\hat{\alpha}\left(|\Gamma(s_i)|{\hat{\varphi}}_{i,k},\eta_{i,l}\right)_{\Lambda_i}
\end{align*}
and into the following vectors:
\begin{align*}
	 & f\in \mathbb{R}^N \text{ s.t. } f_k=\left(f,\varphi_k\right)_{\Omega}\qquad           g_i\in \mathbb{R}^{\hat{N}_i} \text{ s.t. } (g_i)_k=\left(|\Sigma(s_i)|\overline{\overline{g}},\hat{\varphi}_{i,k}\right)_{\Lambda_i}.
\end{align*}
Matrices relative to the various segments $\Lambda_i$, $i=1,\ldots,\I$ are grouped together, forming:
\begin{align*}
	 & \bm{B}=\left[ \bm{B_1}, \bm{B_2},...,\bm{B_{\I}}\right]  \in \mathbb{R}^{N\times N^{\phi}} \qquad \bm{\hat{B}}=\text{diag}\left( \bm{\hat{B}_1},...,\bm{\hat{B}_{\I}}\right)   \in \mathbb{R}^{\hat{N}\times N^{\phi}}                                                         \\
	 & \bm{C}^{\alpha}=\left[ \bm{C_1}^{\alpha}, \bm{C_2}^{\alpha},...,\bm{C_{\I}}^{\alpha}\right] \in \mathbb{R}^{N\times N^{\psi}} \qquad \bm{\hat{C}}^{\alpha}=diag\left(  \bm{\hat{C}_1}^{\alpha},...,\bm{\hat{C}_{\I}}^{\alpha}\right)   \in \mathbb{R}^{\hat{N}\times N^{\psi}}
\end{align*}
being $\hat{N}=\sum_{i=1}^{\I}\hat{N}_i$, $N^{\psi}=\sum_{i=1}^{\I}N_i^{\psi}$ and $N^{\phi}=\sum_{i=1}^{\I}\ N_i^{\phi}$.
Matrices $\bm{\hat{A}_i}$ are grouped as follows, forming matrix $\bm{\hat{A}}$
\begin{displaymath}
	\bm{\hat{A}}=\left[
		\begin{array}{cc}
			\text{diag}\left(  \bm{\hat{A}_{1}},...,\bm{\hat{A}_{\I}}\right) & \bm{Q}^T \\
			\bm{Q}                                                           & \bm{0}
		\end{array}
		\right]
\end{displaymath}
where matrix $\bm{Q}$ simply equates the DOFs placed at the intersections among segments.
We can thus write
\begin{align}
	 & \bm{A}U-\bm{B}\Phi-\bm{C}^{\alpha}\Psi=f\label{eq1discr}                         \\
	 & \bm{\hat{A}}\hat{U}+\bm{\hat{B}}\Phi-\bm{\hat{C}}^{\alpha}\Psi=g\label{eq2discr}
\end{align}
with
\begin{align*}
	 & \hat{U}=\left[\hat{U}_1^T,...,\hat{U}_{\I}^T \right]^T \in \mathbb{R}^{\hat{N}}; \quad g=[g_1^T,g_2^T,...,g_{\I}^T]^T\in\mathbb{R}^{\hat{N}}           \\
	 & \Phi=\left[\Phi_1^T,...,\Phi_{\I}^T \right]^T \in \mathbb{R}^{N^{\phi}}; \quad \Psi=\left[\Psi_1^T,...,\Psi_{\I}^T \right]^T\in \mathbb{R}^{N^{\psi}},
\end{align*}
and, finally setting $W=(U,\hat{U})$,
\begin{equation}
	\bm{\mathcal{A}}=\begin{bmatrix}
		\bm{A} & 0            \\
		0      & \bm{\hat{A}}
	\end{bmatrix}, \qquad
	\bm{\mathcal{B}}=\begin{bmatrix}
		\bm{B} \\
		-\bm{\hat{B}}
	\end{bmatrix},\qquad
	\bm{\mathcal{C}}^{\alpha}=\begin{bmatrix}
		\bm{C}^{\alpha} \\
		\bm{\hat{C}}^{\alpha}
	\end{bmatrix}\qquad
	\mathcal{F}=\begin{bmatrix}
		f \\
		g
	\end{bmatrix},
	\label{Adef}
\end{equation}
the discrete constraint equations are written as:
\begin{equation}
	\bm{\mathcal{A}}W-\bm{\mathcal{B}}\Phi+\bm{\mathcal{C}}^{\alpha}\Psi=\mathcal{F}.\label{eqdiscr_compatta}
\end{equation}

Replacing now the definitions of the discrete variables into the cost functional and replacing the norms in the functional with $L^2$ norms, we can collect the integrals of basis functions into the following matrices
\begin{align*}
	 & \bm{G_i} \in \mathbb{R}^{N \times N} \text{ s.t. } (G_i)_{kl}=\left({\varphi_k}_{|_{\Lambda_i}},{\varphi_l}_{|_{\Lambda_i}}\right)_{\Lambda_i}                 \\
	 & \bm{\hat{G}_i} \in \mathbb{R}^{\hat{N}_i \times \hat{N}_i} \text{ s.t. } (\hat{G}_i)_{kl}=\left({\hat{\varphi}}_{i,k},{\hat{\varphi}}_{i,l}\right)_{\Lambda_i} \\
	 & \bm{G_i^{\psi}} \in \mathbb{R}^{N_i^{\psi} \times N_i^{\psi}}  \text{ s.t. } (G_i^{\psi})_{kl}=\left(\eta_{i,k},\eta_{i,l}\right)_{\Lambda_i}                  \\
	 & \bm{C_i} \in \mathbb{R}^{N\times N_i^{\psi}} \text{ s.t. } (C_i)_{kl}=\left({\varphi_k}_{|_{\Lambda_i}},\eta_{i,l}\right)_{\Lambda_i}                          \\
	 & \bm{\hat{C}_i} \in \mathbb{R}^{\hat{N}_i\times N_i^{\psi}} \text{ s.t. } (\hat{C_i})_{kl}=\left({\hat{\varphi}}_{i,k},\eta_{i,l}\right)_{\Lambda_i}
\end{align*}
and
\begin{equation}
	\label{Gdef}
	\bm{G}=\sum_{i=1}^{\I}\bm{G}_i \in \mathbb{R}^{N \times N} \qquad \bm{\hat{G}}=\text{diag}\left( \bm{\hat{G}_1}^T,...,\bm{\hat{G}_{\I}}^T\right)  \in \mathbb{R}^{\hat{N} \times\hat{N}} \qquad \bm{\mathcal{G}}=\begin{bmatrix}
		\bm{G} & 0 \\0 &\bm{\hat{G}}
	\end{bmatrix}
\end{equation}
\begin{equation*}
	\bm{G^{\psi}}=\text{diag}\left(\bm{G_1^{\psi}},...,\bm{G_{\I}^{\psi}} \right) \in \mathbb{R}^{N^{\psi}\times N^{\psi}}
\end{equation*}
\begin{align*}
	\bm{C}=\left[ \bm{C_1}, \bm{C_2},...,\bm{C_{\I}}\right] \in \mathbb{R}^{N\times N^{\psi}} \quad \bm{\hat{C}}=\text{diag}\left(  \bm{\hat{C}_1},...,\bm{\hat{C}_{\I}}\right)   \in \mathbb{R}^{\hat{N}\times N^{\psi}} \quad \bm{\mathcal{C}}=\begin{bmatrix} \bm{C} \\ \bm{\hat{C}} \end{bmatrix},
\end{align*}
thus deriving the discrete version of the functional, denoted by $\tilde{J}$:
\begin{align}
	\tilde{J} & =\cfrac{1}{2}\left( U^T\bm{G}U-U^T\bm{C}\Psi-\Psi^T\bm{C}^TU+\hat{U}^T\bm{\hat{G}}\hat{U}-\hat{U}^T\bm{\hat{C}}\Psi-\Psi^T\bm{\hat{C}}^T\hat{U}+2\Psi^T\bm{G^{\psi}}\Psi\right)=\nonumber \\
	          & =\cfrac{1}{2}\left( W^T\bm{\mathcal{G}}W-W^T\bm{\mathcal{C}}\Psi-\Psi^T\bm{\mathcal{C}}^TW+2\Psi^T\bm{G^{\psi}}\Psi\right). \label{Jtilde}
\end{align}

The discrete formulation of problem \eqref{minJ} thus is:
\begin{align}
	\min_{(\Phi,\Psi)}\tilde{J}(\Phi,\Psi) \text{ subject to } (\ref{eqdiscr_compatta}). \label{minJtilde}
\end{align}

\section{Resolution method} \label{sec:CG}
\begin{algorithm}
	\SetAlgoLined
	Guess $\mathcal{X}_0=[\Phi_0^T, \Psi_0^T]^T$ \medskip\\
	$r_0=\bm{M}\mathcal{X}_0+d$;\medskip\\
	set $\delta \mathcal{X}_0=-r_0$ and $k=0$;\\
	\While{$\frac{||r_k||}{||d||} >toll$}{
	$\zeta_k=\cfrac{r_k^Tr_k}{\delta \mathcal{X}_{k}^T\bm{M}\delta \mathcal{X}_{k}}$;\\
	$\mathcal{X}_{k+1}=\mathcal{X}_k+\zeta_k\delta \mathcal{X}_{k}$;\\
	$r_{k+1}=r_k+\zeta_k\bm{M}\delta \mathcal{X}_{k}$;\\
	$\beta_{k+1}=\cfrac{r_{k+1}^Tr_{k+1}}{r_{k}^Tr_{k}}$;\\
	$\delta \mathcal{X}_{k+1}=-r_{k+1}+\beta_{k+1}\delta \mathcal{X}_k$;\\
	$k=k+1$;}
	\caption{Conjugate gradient method for $\bm{M}\mathcal{X}+d=0$}
	\label{grad_congAlgo}
\end{algorithm}
The resolution of the previous problem can be efficiently performed via a gradient based method.
By formally replacing
$W=\bm{\mathcal{A}}^{-1}(\bm{\mathcal{B}}\Phi-\bm{\mathcal{C}}^{\alpha}\Psi+\mathcal{F})$
in the functional \eqref{Jtilde}, we obtain
\begin{align*}
	J^\star(\Phi,\Psi) & =\cfrac{1}{2}\Big( (\bm{\mathcal{A}}^{-1}\bm{\mathcal{B}}\Phi+\bm{\mathcal{A}}^{-1}\bm{\mathcal{C}}^{\alpha}\Psi+\bm{\mathcal{A}}^{-1}\mathcal{F})^T\bm{\mathcal{G}}(\bm{\mathcal{A}}^{-1}\bm{\mathcal{B}}\Phi+\bm{\mathcal{A}}^{-1}\bm{\mathcal{C}}^{\alpha}\Psi+\bm{\mathcal{A}}^{-1}\mathcal{F})+\nonumber \\
	                   & \qquad-(\bm{\mathcal{A}}^{-1}\bm{\mathcal{B}}\Phi+\bm{\mathcal{A}}^{-1}\bm{\mathcal{C}}^{\alpha}\Psi+\bm{\mathcal{A}}^{-1}\mathcal{F})^T\bm{\mathcal{C}}\Psi+\nonumber                                                                                                                                        \\&\qquad-\Psi^T\bm{\mathcal{C}}^T(\bm{\mathcal{A}}^{-1}\bm{\mathcal{B}}\Phi+\bm{\mathcal{A}}^{-1}\bm{\mathcal{C}}^{\alpha}\Psi+\bm{\mathcal{A}}^{-1}\mathcal{F})\Big)=\nonumber\\
	                   & =\cfrac{1}{2}~[\Phi^T \quad \Psi^T]\begin{bmatrix}
		\bm{\mathcal{B}}^T\bm{\mathcal{A}}^{-T}\bm{\mathcal{G}}\bm{\mathcal{A}}^{-1}\bm{\mathcal{B}} & \quad \Large\substack{\bm{\mathcal{B}}^T\bm{\mathcal{A}}^{-T}\bm{\mathcal{G}}\bm{\mathcal{A}}^{-1}\bm{\mathcal{C}}^{\alpha}+ \\-\bm{\mathcal{B}}^T\bm{\mathcal{A}}^{-T}\bm{\mathcal{C}}}\\\\
		\Large\substack{(\bm{\mathcal{C}}^{\alpha})^T\bm{\mathcal{A}}^{-T}\bm{\mathcal{G}}\bm{\mathcal{A}}^{-1}\bm{\mathcal{B}}+                                                                                                    \\-\bm{\mathcal{C}}^T\bm{\mathcal{A}}^{-1}\bm{\mathcal{B}}} & \Large\substack{(\bm{\mathcal{C}}^{\alpha})^T\bm{\mathcal{A}}^{-T}\bm{\mathcal{G}}\bm{\mathcal{A}}^{-1}\bm{\mathcal{C}}^{\alpha}+\nonumber \\-\bm{\mathcal{C}}^T\bm{\mathcal{A}}^{-T}\bm{\mathcal{C}}^{\alpha}+\\-(\bm{\mathcal{C}}^{\alpha})^T\bm{A}^{-1}\bm{\mathcal{C}}+2\bm{G^{\psi}}}
	\end{bmatrix}\begin{bmatrix}
		\Phi\smallskip \\\\\\\Psi
	\end{bmatrix}+                                                                                                                                                                                                                      \\[0.5em]
	                   & \qquad+\mathcal{F}^T\begin{bmatrix}
		\bm{\mathcal{A}}^{-T}\bm{\mathcal{G}}\bm{\mathcal{A}}^{-1}\bm{\mathcal{B}} & \quad \bm{\mathcal{A}}^{-T}\bm{\mathcal{G}}\bm{\mathcal{A}}^{-1}\bm{\mathcal{C}}^{\alpha}-\bm{\mathcal{A}}^{-T}\bm{\mathcal{C}}
	\end{bmatrix}\begin{bmatrix}
		\Phi \\\Psi
	\end{bmatrix}+\nonumber                                                                                                                                                                                                                            \\
	                   & \qquad +\cfrac{1}{2}\left( \mathcal{F}^T\bm{\mathcal{A}}^{-T}\bm{\mathcal{G}}\bm{\mathcal{A}}^{-1}\mathcal{F}\right).
\end{align*}
If we set $\mathcal{X}=[\Phi^T,\Psi^T]^T$, we can rewrite $J^\star$ in a compact form as
\begin{equation}
	J^\star(\mathcal{X})=\cfrac{1}{2}\left( \mathcal{X}^T\bm{M}\mathcal{X}+2d^T\mathcal{X}+q\right),\label{Jcompact}
\end{equation}
with
\begin{equation}
	\bm{M}=\begin{bmatrix}
		\bm{\mathcal{B}}^T\bm{\mathcal{A}}^{-T}\bm{\mathcal{G}}\bm{\mathcal{A}}^{-1}\bm{\mathcal{B}} & \quad \bm{\mathcal{B}}^T\bm{\mathcal{A}}^{-T}\bm{\mathcal{G}}\bm{\mathcal{A}}^{-1}\bm{\mathcal{C}}^{\alpha}-\bm{\mathcal{B}}^T\bm{\mathcal{A}}^{-T}\bm{\mathcal{C}} \\\\
		\Large\substack{(\bm{\mathcal{C}}^{\alpha})^T\bm{\mathcal{A}}^{-T}\bm{\mathcal{G}}\bm{\mathcal{A}}^{-1}\bm{\mathcal{B}}+                                                                                                                                           \\-\bm{\mathcal{C}}^T\bm{\mathcal{A}}^{-1}\bm{\mathcal{B}}} & \Large\substack{(\bm{\mathcal{C}}^{\alpha})^T\bm{\mathcal{A}}^{-T}\bm{\mathcal{G}}\bm{\mathcal{A}}^{-1}\bm{\mathcal{C}}^{\alpha}-\bm{\mathcal{C}}^T\bm{\mathcal{A}}^{-T}\bm{\mathcal{C}}^{\alpha}+\\-(\bm{\mathcal{C}}^{\alpha})^T\bm{A}^{-1}\bm{\mathcal{C}}+2\bm{G^{\psi}}}
	\end{bmatrix} \label{CGmatrix}
\end{equation}\smallskip\\
\begin{equation}
	d^T=
	\mathcal{F}^T\begin{bmatrix}
		\bm{\mathcal{A}}^{-T}\bm{\mathcal{G}}\bm{\mathcal{A}}^{-1}\bm{\mathcal{B}} & \quad \bm{\mathcal{A}}^{-T}\bm{\mathcal{G}}\bm{\mathcal{A}}^{-1}\bm{\mathcal{C}}^{\alpha}-\bm{\mathcal{A}}^{-T}\bm{\mathcal{C}}
	\end{bmatrix},
\end{equation}
\begin{equation}
	q=\mathcal{F}^T\bm{\mathcal{A}}^{-T}\bm{\mathcal{G}}\bm{\mathcal{A}}^{-1}\mathcal{F}.
\end{equation}
Matrix $\bm{M}$ is symmetric positive definite, as it follows from the structure of functional \eqref{Jtilde} and from the equivalence of this formulation with the well posed problem \eqref{minJtilde} \cite{3D1Darxiv}.
The minimum of \eqref{Jcompact} is given by condition
\begin{equation}
	\nabla J^\star=\bm{M}\mathcal{X}+d=0. \label{grad_cong_system}
\end{equation}The minimization of the unconstrained problem \eqref{Jcompact} can be performed via a conjugate gradient method, as reported in Algorithm~\ref{grad_congAlgo}.

Let us observe that the application of matrix $\bm{M}$ to an array, say $\delta\mathcal{X}$, does not involve the explicit computation of matrix $\bm{M}$ and of the inverse matrix $\bm{\mathcal{A}}^{-1}$ . Indeed the quantity $\bm{M}\delta\mathcal{X}$, whose computation is required several times in Algorithm~\ref{grad_congAlgo}, can be performed as:
\begin{equation*}
	\bm{M}\delta\mathcal{X}=\begin{bmatrix}
		\bm{\mathcal{B}}^T\delta \mathcal{P} \\ (\bm{\mathcal{C}}^{\alpha})^T\delta \mathcal{P}-\bm{\mathcal{C}}^{T}\delta W+2\bm{G^{\psi}}\delta \Psi
	\end{bmatrix}
\end{equation*}
where $\delta \mathcal{P}$ is obtained as the solution of the system
\begin{displaymath}
	\bm{\mathcal{A}}^{T}\delta \mathcal{P}=\bm{\mathcal{G}}\delta W -\bm{\mathcal{C}}\delta\Psi,
\end{displaymath}
which, in virtue of the structure of matrix $\bm{\mathcal{A}}$, requires the resolution of independent sub-problems on each of the 1D segments and on the 3D domain.

\section{Numerical results}\label{Num_res}

In this section we propose three numerical tests to show the applicability and the performances of the proposed conjugate gradient solver for the optimization formulation of coupled 3D-1D problems. The first test proposes a comparison between the solution of a fully 3D-3D simulation on a conforming mesh and the solution of the corresponding reduced 3D-1D problem with the proposed approach. The quality of the solution is evaluated in terms of total flux conservation. The second test takes into account the problem of the computation of the equivalent permeability of a porous medium when crossed by a set of conductive small channels. Finally, the third test shows the potential of the approach in dealing with extremely complex configurations, considering a set of $1000$ possibly intersecting segments embedded in a porous matrix.

Simulations are performed using linear Lagrangian finite elements on tetrahedral meshes for the 3D domain, whereas linear Lagrangian finite elements on equally spaced meshes are used on each segment $\Lambda_i$ $i=1,...\I$ for the unknowns $\hat{U}$ and $\Psi$.
Piecewise constant basis functions on equally spaced nodes are instead used for $\Phi$ on each segment. For simplicity, mesh refinement is denoted by means of a unique parameter $h$, representing the maximum diameter of the tetrahedra for the 3D mesh of $\Omega$. The refinement level of the 1D meshes is related to $h$ as follows: called $N_i^\star$ the number of intersection points between the faces of the tetrahedra of the 3D mesh and segment $\Lambda_i$, we build on $\Lambda_i$ a mesh made of $N_i^\star$ equally spaced nodes for variable $\hat{U}$ and $\frac12 N_i^\star$ equally spaced nodes for variables $\Psi$ and $\Phi$. Clearly different refinement levels could be chosen on each segment and for each 1D unknown. Such analysis, however, is out of the scope of the present work; the interested reader can refer to \cite{3D1Darxiv} for more detail on this issue. Parameters $\alpha$ and $\hat{\alpha}$ are set to one for all the simulations, even if any other strictly positive value can be used.

\subsection{Problem 1: comparison with a 3D-3D simulation}
\begin{figure}
	\centering
	\begin{subfigure}{0.48\textwidth}
		\includegraphics[width=\textwidth]{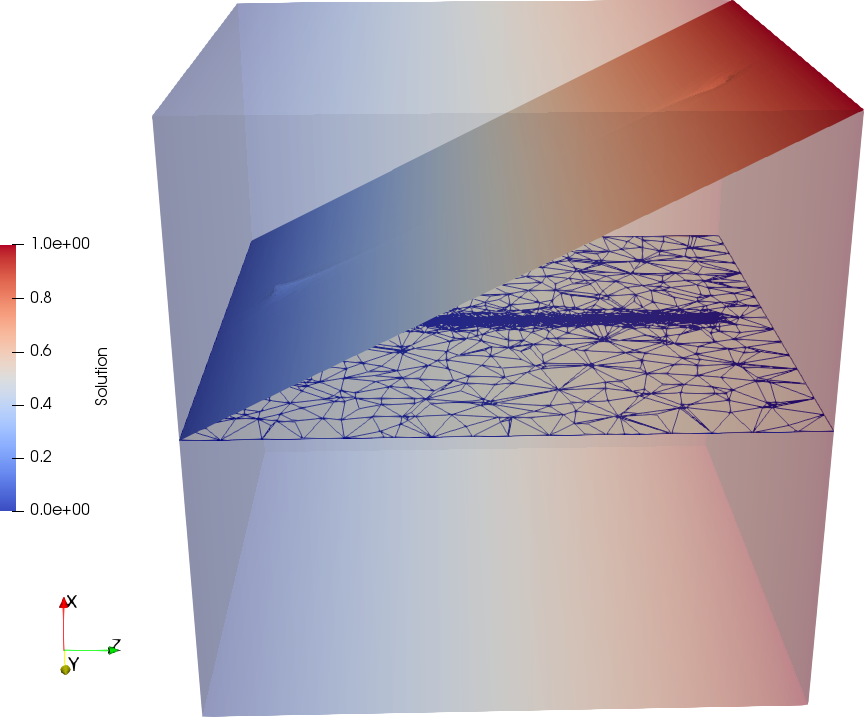}
		\caption{Equi-dimensional problem}
		\label{3D3DDomainInsideSeg}
	\end{subfigure}
	\begin{subfigure}{0.48\textwidth}
		\includegraphics[width=\textwidth]{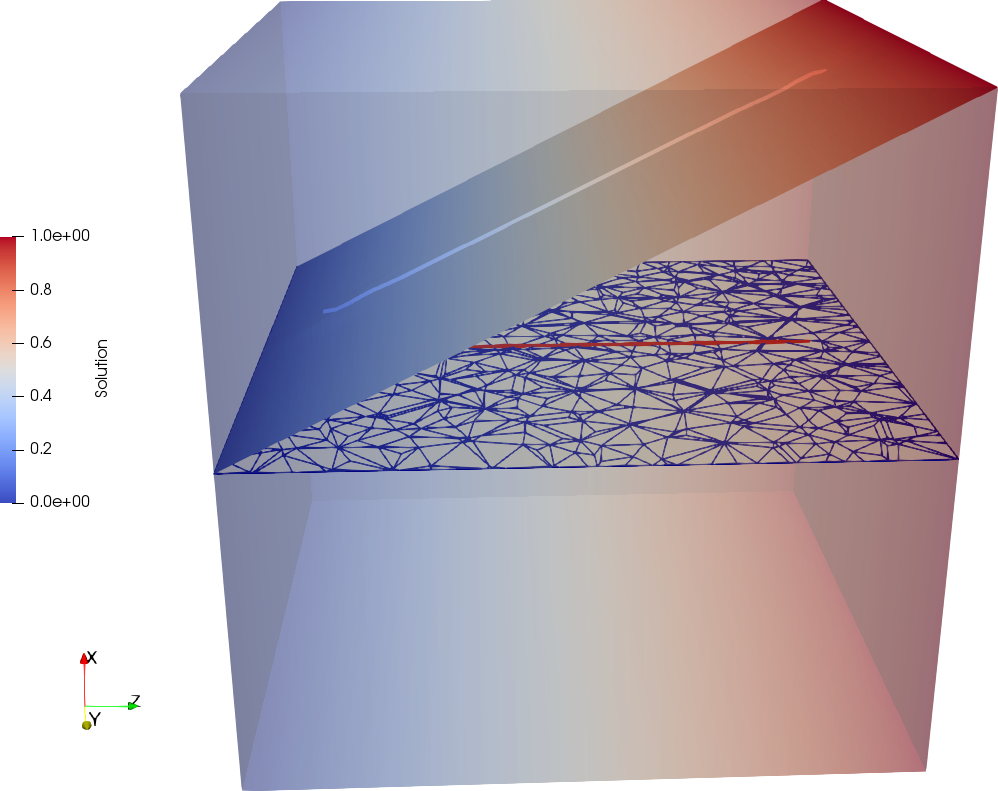}
		\caption{3D-1D reduced problem}
		\label{3D1DDomainInsideSeg}
	\end{subfigure}
	\caption{Problem 1 - Coarse mesh and solution on a plane containing the centreline of the inclusion}
	\label{fig:1:result}
\end{figure}
\begin{table}
	\caption{Problem 1 - Comparison with the 3D-3D case: number of DOFs and fluxes across the faces}
	\label{table_3D3D}
	\renewcommand*{\arraystretch}{1.2}
	\begin{center}
		\begin{tabular}{ccc|cc}
			\hline
			                          & \multicolumn{2}{c|}{$\bm{h=1.0\cdot10^{-1}}$} & \multicolumn{2}{c}{$\bm{h=4.6\cdot10^{-2}}$}                                              \\
			\hline
			                          & \textbf{3D-1D}                                & \textbf{3D-3D}                               & \textbf{3D-1D}       & \textbf{3D-3D}      \\
			\hline
			$\bm{N}$                  & $2998$                                        & $11295$                                      & $26109$              & $36343$             \\
			$\bm{\hat{N}}$            & $37$                                          & $-$                                          & $87$                 & $-$                 \\
			\hline
			$\bm{\sigma_1}$ (outflow) & $2.0117$                                      & $2.0112$                                     & $2.0116$             & $2.0108$            \\
			$\bm{\sigma_2}$           & $8.60\cdot 10^{-7}$                           & $-2.03\cdot 10^{-6}$                         & $-1.06\cdot10^{-6}$  & $1.15\cdot10^{-6}$  \\
			$\bm{\sigma_3}$           & $-1.68\cdot 10^{-5}$                          & $4.80\cdot 10^{-6}$                          & $1.10\cdot10^{-7}$   & $7.28\cdot 10^{-7}$ \\
			$\bm{\sigma_4}$           & $2.81\cdot 10^{-6}$                           & $-2.27\cdot 10^{-6}$                         & $-4.21\cdot10^{-7}$  & $-2.81\cdot10^{-6}$ \\
			$\bm{\sigma_5}$           & $-4.14\cdot 10^{-6}$                          & $-6.97\cdot10^{-6}$                          & $-2.78\cdot 10^{-6}$ & $-1.26\cdot10^{-6}$ \\
			$\bm{\sigma_6}$ (inflow)  & $-2.0120$                                     & $-2.0107$                                    & $-2.0116$            & $-2.0109$           \\
			\hline
		\end{tabular}
	\end{center}
\end{table}
The first example takes into account a simple setting, with a single inclusion lying in the interior of a cubic domain. A comparison is proposed between the solution obtained solving the equi-dimensional 3D-3D problem via a conforming mesh, and the solution of the reduced 3D-1D problem on a non-conforming mesh via the proposed approach.
Let us consider a cube of edge $l=2$ whose barycentre is located at the origin of a reference system $xyz$, and a segment $\Lambda$ lying on the $z$-axis and going from $z=-0.8$ to $z=0.8$. This segment is supposed to be the centreline of a cylindrical channel of radius $\check{R}=10^{-2}$ and transmissivity $\bm{\tilde{K}}=10^2$, while in the cube we consider a permeability coefficient $\bm{K}=1$. Let us impose homogeneous Neumann conditions on all the lateral faces of the cube, and Dirichlet boundary conditions on the top and bottom faces, respectively equal to 1 and 0. Homogeneous Neumann conditions are also imposed at segment endpoints lying in the interior of the domain.

In the equi-dimensional setting, the cylindrical inclusion is approximated by a prism with $16$ faces and the mesh is conforming at the interface between the inclusion and the outer domain. The resulting mesh is thus refined towards the inclusion, in order to match the edge-size of the elements on the interfaces as shown in Figure~\ref{3D3DDomainInsideSeg}, where such adapted mesh is shown by its intersection with the plane containing the centreline of the inclusion and normal to the $x$-axis.
For the 3D-1D problem, instead, the inclusion is reduced to its centreline, which arbitrarily crosses the elements of the 3D mesh, see Figure~\ref{3D1DDomainInsideSeg}. Figures~\ref{3D3DDomainInsideSeg}-\ref{3D1DDomainInsideSeg} also provide a plot of the solution on the same plane.

Let us denote by $\bm{\sigma_i}=-\int_{\partial\Omega_i}K\nabla U\cdot \bm{n}_i$ the amount of flux leaving the $i$-th face of the 3D domain, being $\bm{n}_i$ the outward pointing normal vector to face $\partial \Omega_i$, $i=1,...,6$. We analyze the performances of our 3D-1D reduced model by comparing the computed fluxes with the ones obtained with the 3D-3D simulation on two different meshes for the 3D domain.
The results are collected in Table~\ref{table_3D3D}. A coarse mesh with $h=1\cdot10^{-1}$ and a fine mesh with $h=4.6\cdot10^{-2}$ are considered. Since the mesh for the equi-dimensional case is adapted at the interface, mesh size close to the inclusion is constrained by the conformity requirement and not by mesh parameter $h$. The number of the degrees of freedom $\bm{N}$ is also provided in Table~\ref{table_3D3D} and can be used to compare the refinement level of the meshes of the different approaches. We can observe that the results carried out by the proposed approach for the reduced problem are in line with the ones obtained by solving the equi-dimensional problem. In particular, the weak approximation of the homogeneous neumann boundary conditions is comparable between the two solutions and also the value of the influx and outflux is in good agreement. Denoting by $\bm{\sigma}^{tot}:=\left| \sum_i \bm{\sigma_i}\right| $ the total flux mismatch, we obtain values of $3.48\cdot10^{-4}$ on the coarse mesh and $3.00\cdot 10^{-5}$ on the fine mesh for the solution of the reduced problem and values of $4.73\cdot10^{-4}$ and $9.00\cdot 10^{-5}$ on the coarse and fine meshes for the solution of the equi-dimensional problem.

\subsection{Problem 2: computation of equivalent transmissivity}
\newcommand{\segq}{\textbf{Seg40 }}
\newcommand{\sego}{\textbf{Seg80 }}
\newcommand{\Keq}{\bm{K_{\textit{eq}}}}
\newcommand{\out}{\text{out}}
\newcommand{\inn}{\text{in}}
\begin{figure}
	\centering
	\begin{subfigure}{0.44\textwidth}
		\includegraphics[width=\textwidth]{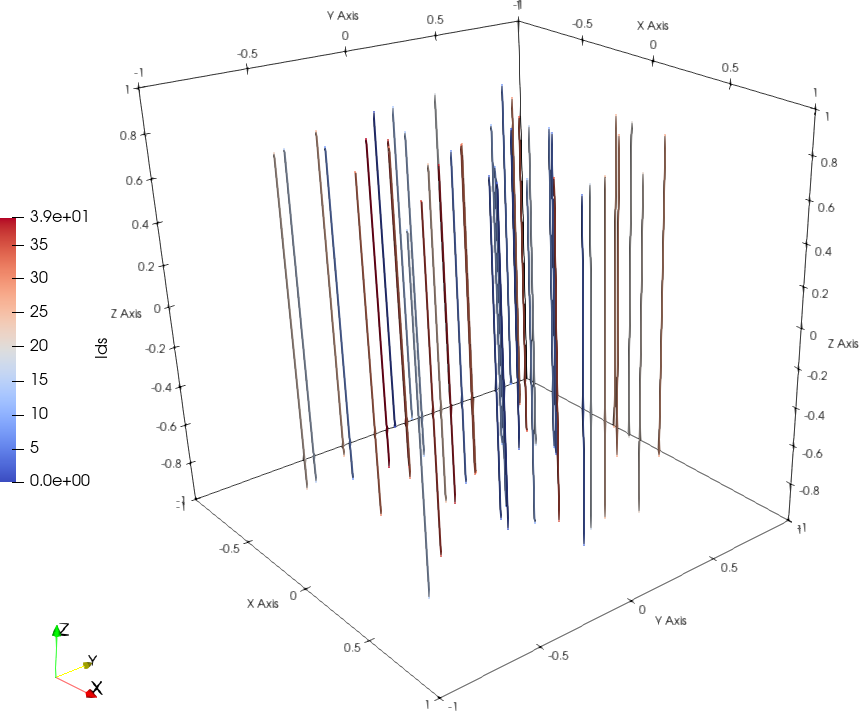}
		\caption{\segq configuration}
		\label{fig:40:domain}
	\end{subfigure}
	\begin{subfigure}{0.44\textwidth}
		\includegraphics[width=\textwidth]{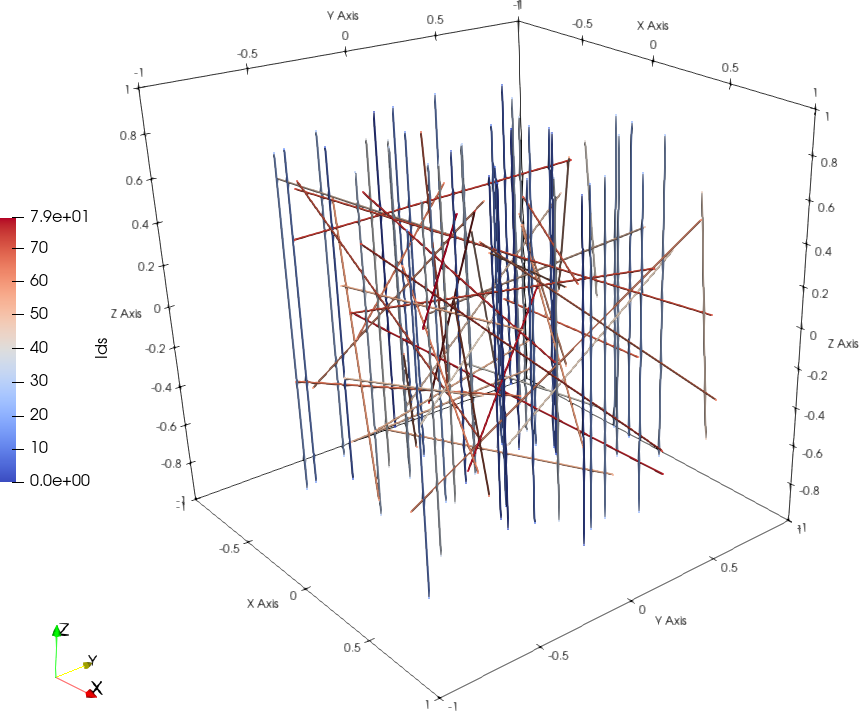}
		\caption{\sego configuration}
		\label{fig:80:domain}
	\end{subfigure}
	\caption{Problem 2 - \segq and \sego configurations}
	\label{fig:40-80:domain}
\end{figure}
\begin{figure}
	\centering
	\includegraphics[width=0.5\textwidth]{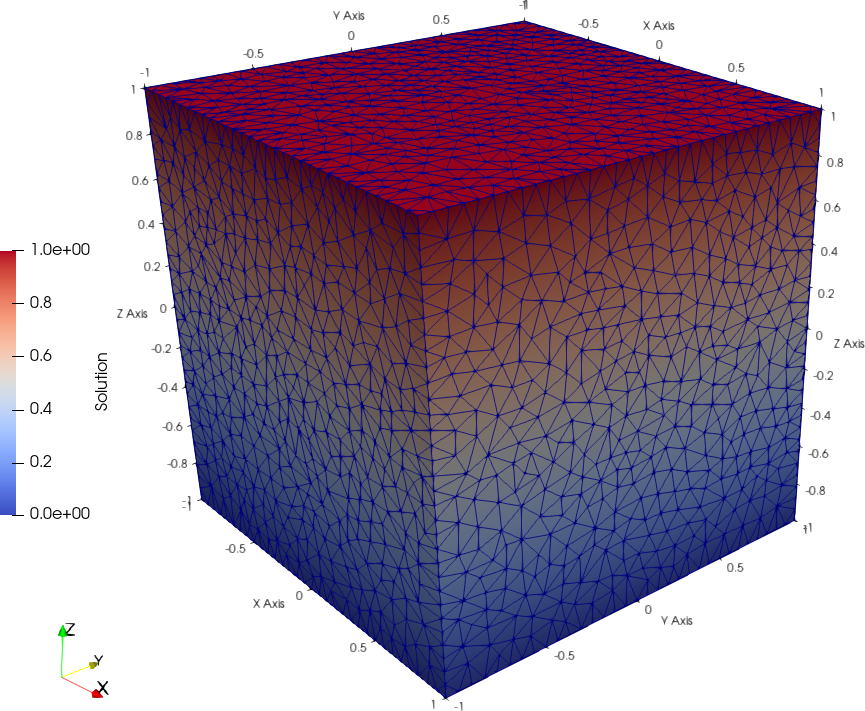}
	\caption{Problem 2 - Solution on the full domain for the \sego setting}
	\label{fig:80:solution}
\end{figure}
\begin{figure}
	\centering
	\begin{subfigure}{0.44\textwidth}
		\includegraphics[width=\textwidth]{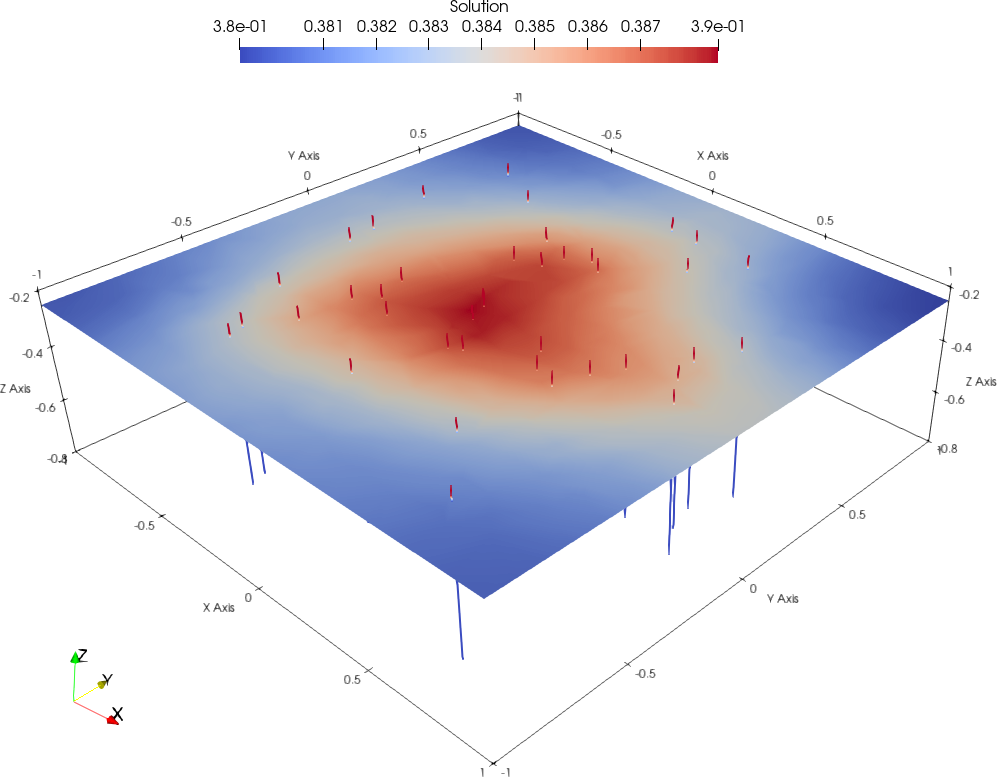}
		\caption{\segq}
		\label{fig:80:solution:block}
	\end{subfigure}
	\begin{subfigure}{0.44\textwidth}
		\includegraphics[width=\textwidth]{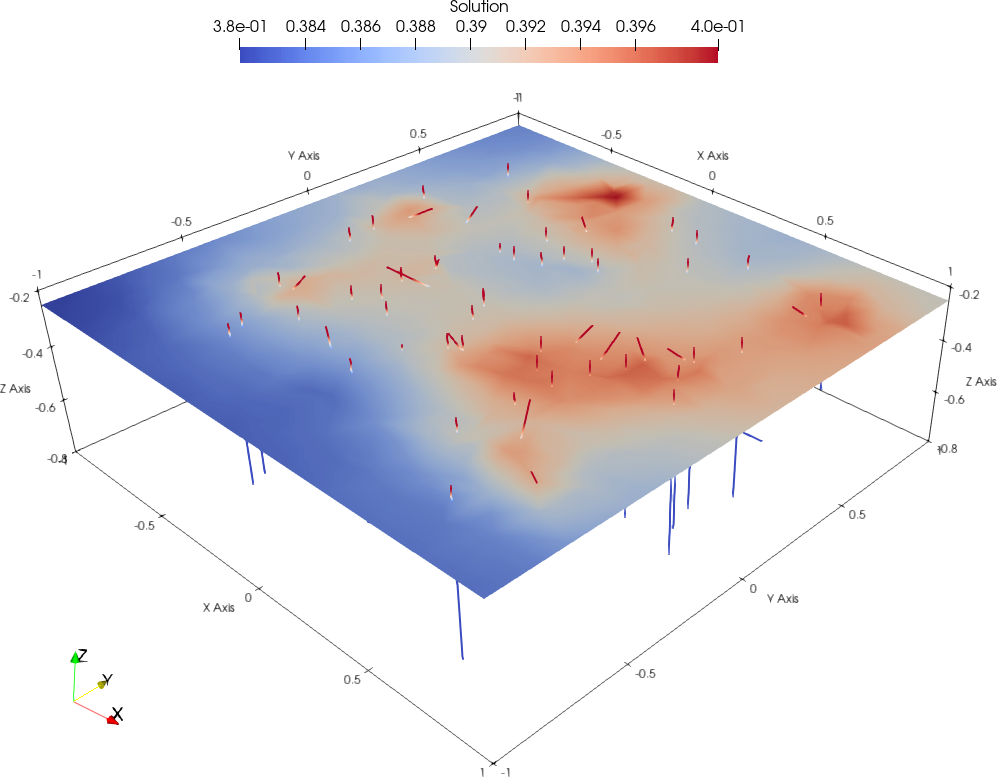}
		\caption{\sego}
		\label{fig:80:solution:cut}
	\end{subfigure}
	\caption{Problem 2 - Solution on a section of the domain with a plane orthogonal to the $z$-axis located at $z = -0.25$}
	\label{fig:40-80:solution:cut}
\end{figure}
\begin{table}
	\caption{Problem 2 - Values of outlet flux, relative flux mismatch and $\Keq$ for the two considered settings}
	\label{Flux_table}
	\renewcommand*{\arraystretch}{1.4}
	\begin{center}
		\begin{tabular}{ccccc}
			\hline
			                              & \textbf{flux dir.}     & $\bm{|\sigma^{\out}|}$ & $\cfrac{\bm{||\sigma^{\out}|-|\sigma^{\inn}||}}{\bm{|\sigma^{\out}|}}$ & $\bm{K_{eq}}$   \\
			\hline
			\multirow{2}{*}{$\bm{Seg40}$} & $\bm{z}$\textbf{-axis} & $2.44$                 & $2.11 \cdot 10^{-5}$                                                   & $K_{eq}^z=1.22$ \\
			                              & $\bm{x}$\textbf{-axis} & $2.00$                 & $3.64\cdot 10^{-7}$                                                    & $K_{eq}^x=1.00$ \\
			\hline
			\multirow{2}{*}{$\bm{Seg80}$} & $\bm{z}$\textbf{-axis} & $2.55$                 & $1.60\cdot 10^{-4}$                                                    & $K_{eq}^z=1.28$ \\
			                              & $\bm{x}$\textbf{-axis} & $2.10$                 & $4.01\cdot 10^{-4}$                                                    & $K_{eq}^x=1.05$ \\
			\hline
		\end{tabular}
	\end{center}
\end{table}
Let us consider the same cube of edge $l=2$ that was introduced for the previous numerical example, and a set of $\I$ segments $\left\lbrace \Lambda_i \right\rbrace_{i=1}^\I$. In a first configuration, labeled \segq, we have $\I=40$ and all the segments are parallel to the $z$-axis and go from $z=-0.8$ to $z=0.8$. The location on the $xy$-plane is randomly generated from a uniform distribution, with $-0.8<x,y<0.8$ (see Figure~\ref{fig:40-80:domain} on the left). As in the previous numerical example, we suppose these segments to be the reduction to the centreline of 40 cylinders $\Sigma_i$ of radius $\check{R}_i=10^{-2}$ and transmissivity $\bm{\tilde{K}}=10^2$, whereas in the cube we consider again a permeability coefficient $\bm{K}=1$.
A second configuration is also considered, called \sego, in which 40 additional segments with random orientation and position in space are added to the \segq setting. Even these segments are supposed to be the reduction to the centreline of cylinders of radius  $10^{-2}$ and transmissivity $\bm{\tilde{K}}=10^2$. Their extremes are contained in a box with $-0.8<x,y,z<0.8$ (see Figure~\ref{fig:40-80:domain}).
We compute the equivalent transmissivity $\Keq$ of an homogenized material, resulting from the presence of the inclusions, by the proposed gradient based scheme for the optimization approach. We expect this material to be anisotropic as, for both settings, at least 40 segments are all oriented in the same direction: for this reason we compare the equivalent transmissivity in the $z$ direction and the one along an orthogonal direction, namely the $x$-direction, denoting them by $K_{eq}^z$ and $K_{eq}^x$, respectively. In order to compute $\Keq^z$ we impose Dirichlet boundary conditions on the top and on the bottom faces of the cube, prescribing a unitary pressure drop, whereas we consider homogeneous Neumann conditions on the other faces. This means that the top face will be the flux inlet face, while the bottom face will be the outlet. To compute $\Keq^x$ we impose, instead a unitary pressure drop between the two faces of the cube orthogonal to the $x$-axis, with the inlet face at $x=1$ and the outlet face at $x=-1$, and no flux conditions on the other faces. In both cases we impose homogeneous Neumann conditions at all segment endpoints. Let us denote by $\bm{\sigma}^{\out}$ the flux leaving the cube from the outlet face $\partial \Omega_{\out}$, of area $|\partial \Omega_{\out}|$ and outward unit normal vector $\bm{n}^{\out}$. We thus have
\begin{equation}
	\Keq=\cfrac{|\sigma^{\out}|}{|\partial\Omega^{\out}|\cdot0.5}\label{Keqz}
\end{equation}
with $\bm{\sigma}^{\out}=-\int_{\partial\Omega^{\out}}K\nabla U\cdot \bm{n}^{\out}$, being $0.5$ the value of the average pressure gradient across the cube in the flux direction.

As an example, the solution obtained for the \sego setting on a mesh with parameter $h=0.086$ is shown in Figure~\ref{fig:80:solution}, whereas Figure~\ref{fig:40-80:solution:cut} shows a section of the solution on a plane orthogonal to the $z$-axis located at $z = -0.25$, for both settings, on the same mesh. We can see how the inclusions alter the pattern of the solution.
The obtained results are collected in Table~\ref{Flux_table} which, in particular, reports the amount of flux $|\sigma^{\out}|$ leaving the cube from the outlet face, the relative mismatch between $\sigma^{\out}$ and the flux $\sigma^{\inn}=-\int_{\partial\Omega^{\inn}}\bm{K}\nabla U\cdot \bm{n}^{\inn}$ entering from the inlet face, and the computed values of $\Keq^z$ and $\Keq^x$. We remark that the relative mismatch $\cfrac{\bm{||\sigma^{\out}|-|\sigma^{\inn}||}}{\bm{|\sigma^{\out}|}}$ can be used as a proxy for solution accuracy. Data is obtained for a mesh parameter $h=0.086$.
As expected, the presence of a set of parallel vessels along the flux direction leads to an equivalent transmissivity $\Keq^z$ higher than the permeability of the porous medium alone. On the contrary, the value of $\Keq^x$ remains equal to $\bm{K}$ for the \segq setting, as expected given the orientation of the inclusions, whereas it is slightly increased by the presence of the additional segments with random orientation in the \sego configuration. We can observe that, in all cases, very small values of relative flux mismatch are observed, in line with the values obtained for Problem 1.

\subsection{Problem 3: multiple inclusions - 1000 segments}
\begin{figure}
	\centering
	\includegraphics[width=0.5\textwidth]{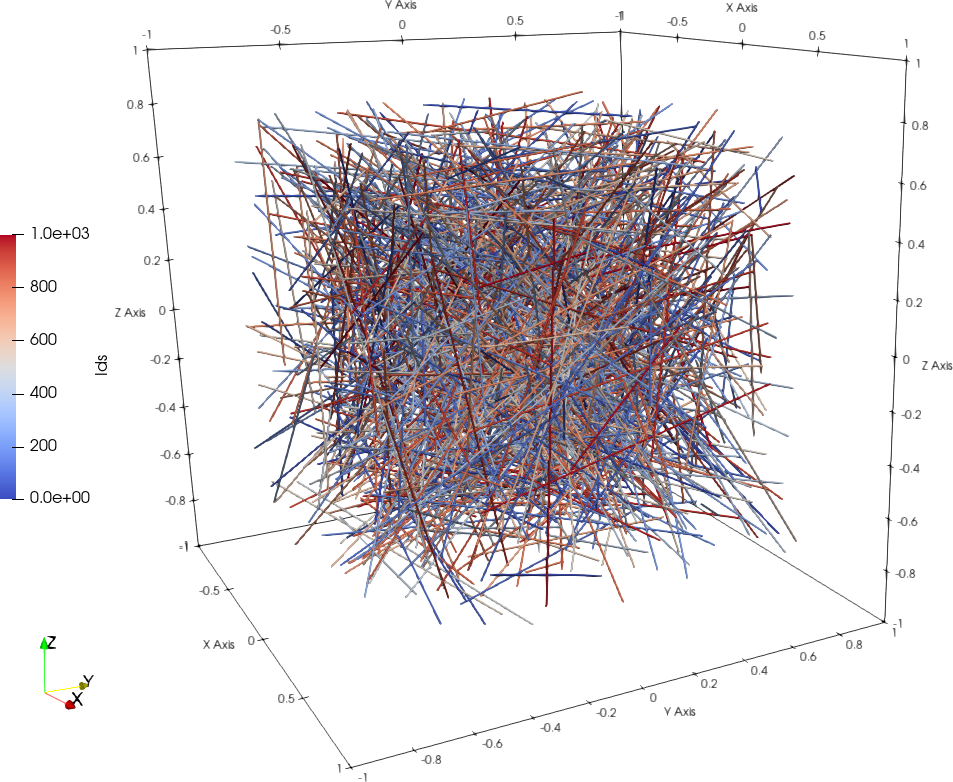}
	\caption{Problem 3 - domain segments}
	\label{fig:1000:domain}
\end{figure}
\begin{figure}
	\centering
	\begin{subfigure}{0.40\textwidth}
		\includegraphics[width=\textwidth]{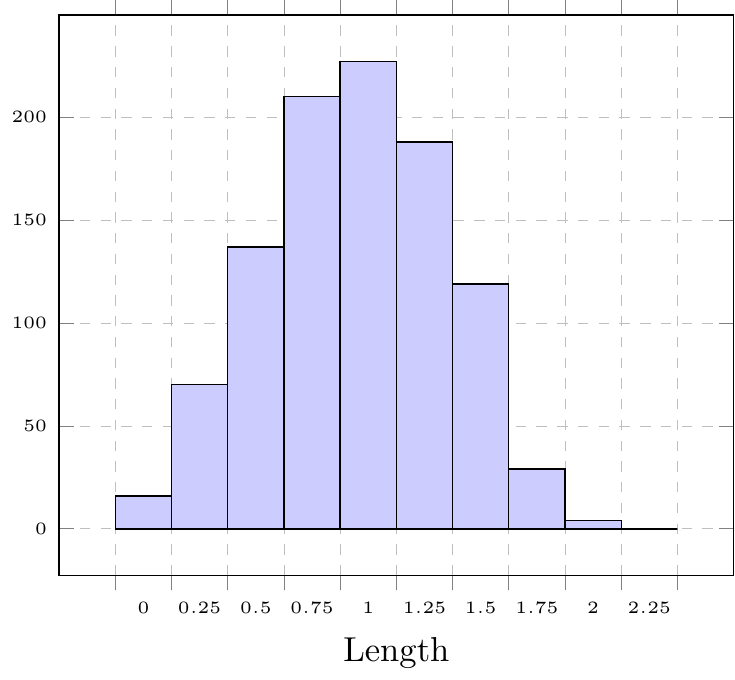}
		\caption{Length distribution}
		\label{fig:1000:statistics:length}
	\end{subfigure}
	\begin{subfigure}{0.40\textwidth}
		\includegraphics[width=\textwidth]{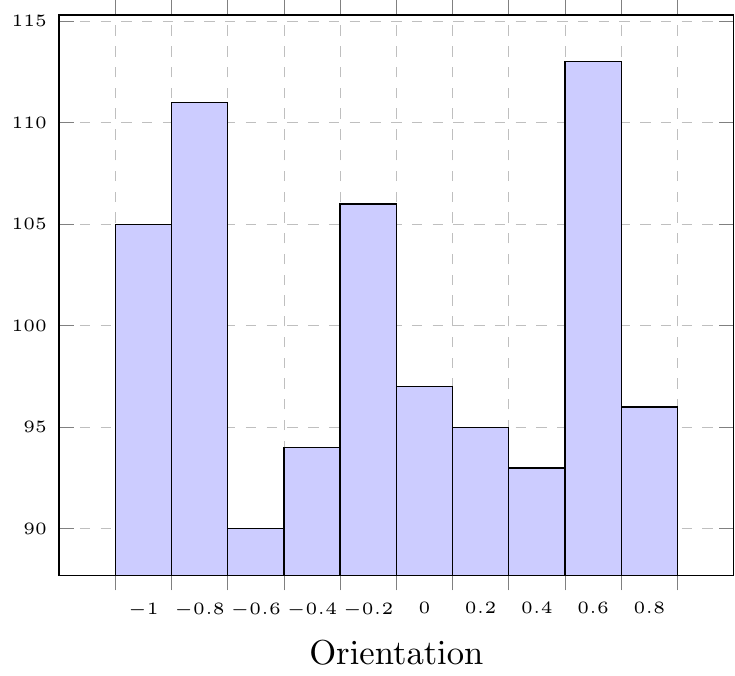}
		\caption{Orientation distribution measured as segment tangent vector $t_s$ dot z-axis unit vector $k$}
		\label{fig:1000:statistics:orientation}
	\end{subfigure}
	\caption{Problem 3 - statistics of domain segments}
	\label{fig:1000:statistics}
\end{figure}
\begin{figure}
	\centering
	\begin{subfigure}{0.44\textwidth}
		\includegraphics[width=\textwidth]{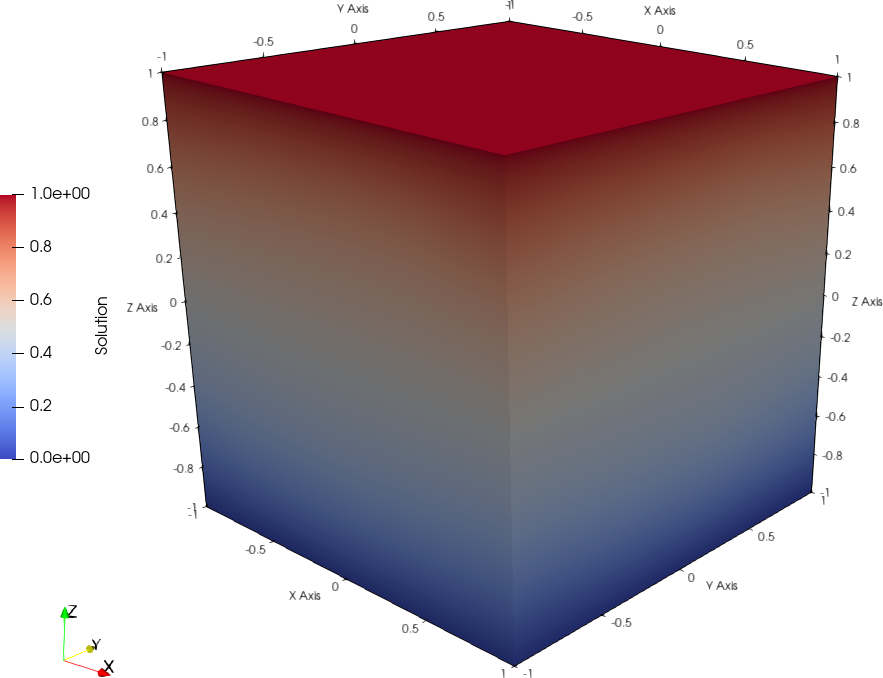}
		\caption{Solution on the whole 3D domain}
		\label{fig:1000:solution:block}
	\end{subfigure}
	\begin{subfigure}{0.44\textwidth}
		\includegraphics[width=\textwidth]{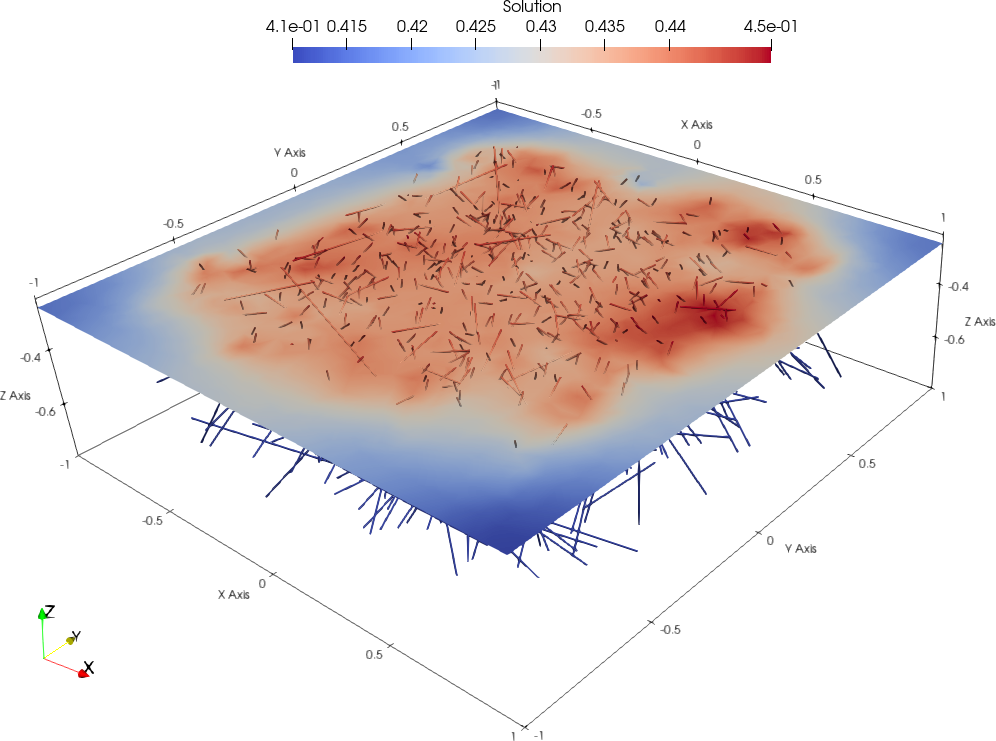}
		\caption{Solution on a section of the domain with a plane orthogonal to the $z$-axis located at $z = -0.25$}
		\label{fig:1000:solution:cut}
	\end{subfigure}
	\caption{Problem 3 - Example solutions on the whole 3D domain and on a section of the domain}
	\label{fig:1000:solution}
\end{figure}
\begin{figure}
	\centering
	\begin{subfigure}{0.32\textwidth}
		\includegraphics[width=\textwidth]{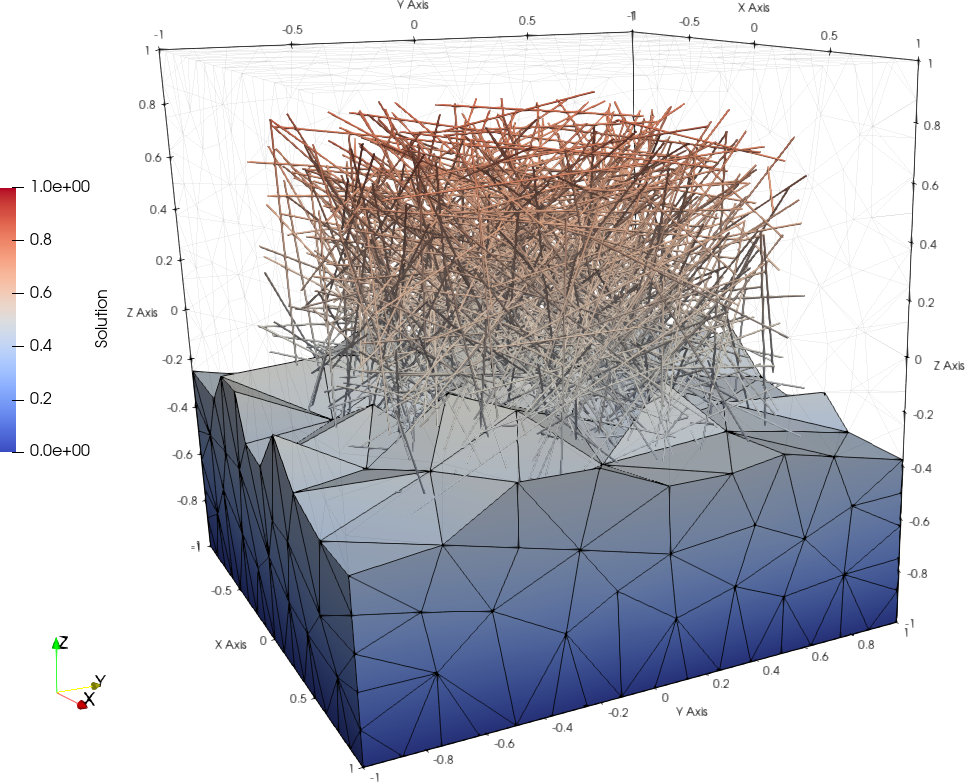}
		\caption{Coarse mesh $h=10^{-\frac{2}{3}}$}
		\label{fig:1000:solution:coarse}
	\end{subfigure}
	\begin{subfigure}{0.32\textwidth}
		\includegraphics[width=\textwidth]{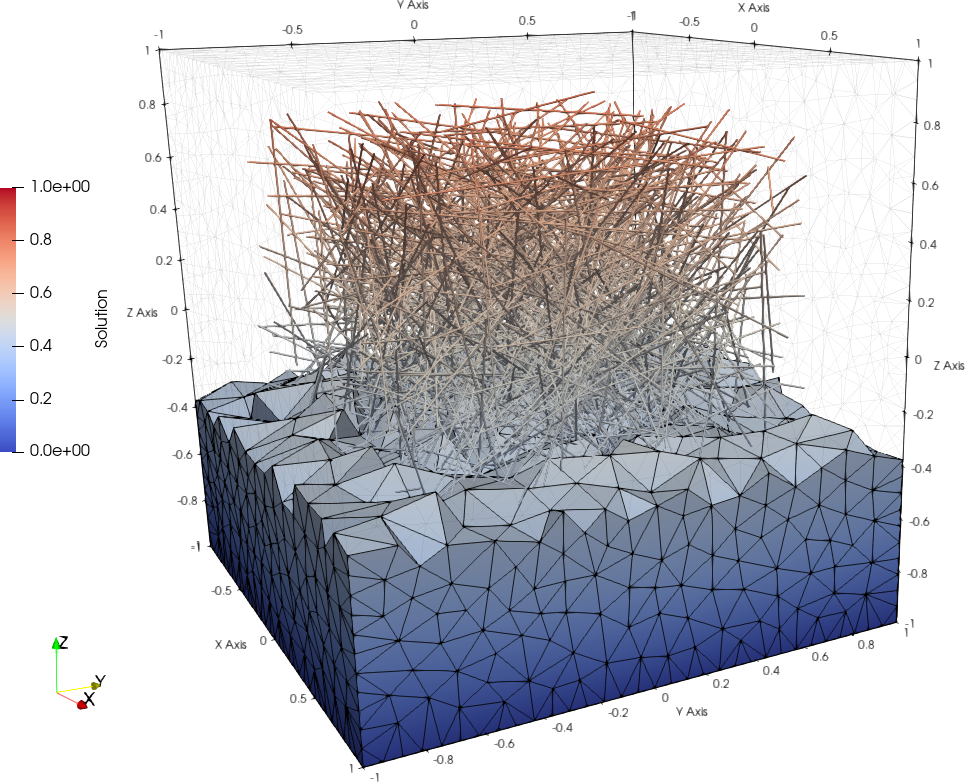}
		\caption{Mean mesh $h=10^{-1}$}
		\label{fig:1000:solution:mean}
	\end{subfigure}
	\begin{subfigure}{0.32\textwidth}
		\includegraphics[width=\textwidth]{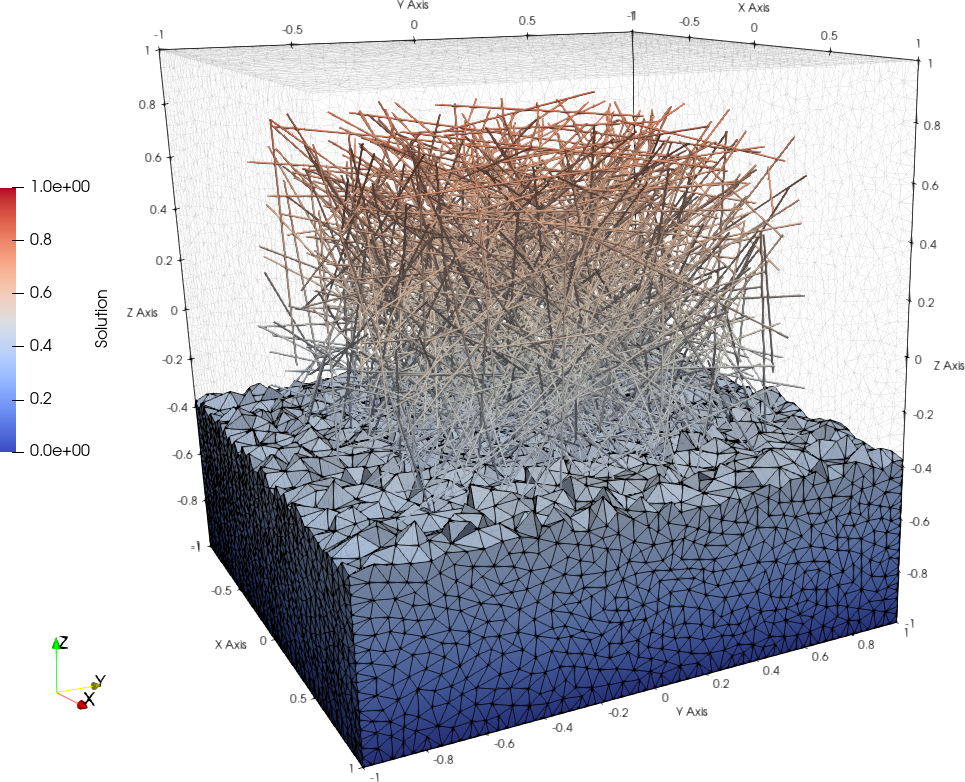}
		\caption{Fine mesh $h=10^{-\frac{4}{3}}$}
		\label{fig:1000:solution:fine}
	\end{subfigure}
	\caption{Problem 3 - details of the three used computational meshes}
	\label{fig:1000:solution:mesh}
\end{figure}
\begin{table}
	\caption{Problem 3 - DOFs and CG iterations}
	\label{fig:1000:dofs}
	\renewcommand*{\arraystretch}{1.4}
	\begin{center}
		\begin{tabular}{ccccccc}
			\hline
			                & $h$                 & $N$     & $\hat{N}$ & $N^{\phi}$ & $N^{\psi}$ & $CG_{It} / (N^{\phi} + N^{\psi})$ \\
			\hline
			\textbf{Coarse} & $10^{-\frac{2}{3}}$ & $400$   & $14885$   & $6695$     & $7695$     & $2.4 \cdot 10^{-1}$               \\
			\hline
			\textbf{Mean}   & $10^{-1}$           & $2998$  & $28351$   & $13426$    & $14426$    & $2.2 \cdot 10^{-1}$               \\
			\hline
			\textbf{Fine}   & $10^{-\frac{4}{3}}$ & $26109$ & $58736$   & $28604$    & $29604$    & $2.4 \cdot 10^{-1}$               \\
			\hline
		\end{tabular}
	\end{center}
\end{table}
\begin{table}
	\caption{Problem 3 - fluxes data}
	\label{fig:1000:fluxes}
	\renewcommand*{\arraystretch}{1.4}
	\begin{center}
		\begin{tabular}{cccc}
			\hline
			                & $\bm{|\sum \sigma^i|}$ & $\bm{|\sigma^{\out}|}$ & $\cfrac{\bm{||\sigma^{\out}|-|\sigma^{\inn}||}}{\bm{|\sigma^{\out}|}}$ \\
			\hline
			\textbf{Coarse} & $1.17 \cdot 10^{-2}$   & $3.2482$               & $1.02 \cdot 10^{-3}$                                                   \\
			\textbf{Mean}   & $3.30 \cdot 10^{-3}$   & $3.6289$               & $4.53 \cdot 10^{-4}$                                                   \\
			\textbf{Fine}   & $9.90 \cdot 10^{-4}$   & $3.6053$               & $6.31 \cdot 10^{-5}$                                                   \\
			\hline
		\end{tabular}
	\end{center}
\end{table}
The last proposed problem, takes into account a set of $1000$ segments embedded in a cubic block of porous material. As in the previous case, the cube has edge length equal to 2 and its barycenter is placed at the origin of a reference system $xyz$. Segments are randomly oriented in the 3D space, as detailed in Figure~\ref{fig:1000:domain} and Figure~\ref{fig:1000:statistics}. A unitary pressure drop is imposed between the top and bottom faces of the domain, all other faces being instead insulated, as well as the extreme of the segments. Simulations are performed on three meshes: a coarse mesh with parameter $h=10^{-\frac{2}{3}}$, an intermediate mesh with $h=10^{-1}$ and a fine mesh with $h=10^{-\frac{4}{3}}$, as shown in Figure~\ref{fig:1000:solution:mesh}. The corresponding numbers of degrees of freedom are reported in Table~\ref{fig:1000:dofs}. The table also reports the number of iterations required by the conjugate gradient scheme, relative to the number of unknowns of the unconstrained problem, to solve the problem up to a relative residual of $10^{-6}$. We can see that the number of iteration is quite stable with respect to mesh refinement. Nonetheless, a preconditioner could be used in order to further reduce the number of iterations, but this is deferred to a forthcoming work. The global flux mismatch is reported in Table~\ref{fig:1000:fluxes} as a proxy of solution accuracy.

\section{Conclusions} \label{sec:Concl}
A gradient based resolution scheme is here proposed for the PDE-constrained optimization approach for coupled 3D-1D problems.
An equivalent unconstrained formulation of the minimization problem is derived and the application of the conjugate gradient scheme to such problem is described and discussed. Numerical examples on quite complex configurations show the applicability and effectiveness of the approach.

\section*{Acknowledgements}
This work is supported by the MIUR project ``Dipartimenti di Eccellenza 2018-2022'' (CUP E11G18000350001), PRIN project ``Virtual Element Methods: Analysis and Applications'' (201744KLJL\_004) and by INdAM-GNCS. Computational resources are partially supported by SmartData@polito.

\bibliographystyle{elsarticle-num}
\bibliography{3D1D_DeniseBib}
\end{document}